# MINIMUM BRAIDS: A COMPLETE INVARIANT OF KNOTS AND LINKS

THOMAS A. GITTINGS
*1221-C Central Street*
*Evanston, IL 60201*
*tomgittings@aol.com*

## ABSTRACT

Minimum braids are a complete invariant of knots and links. This paper defines minimum braids, describes how they can be generated, presents tables for knots up to ten crossings and oriented links up to nine crossings, and uses minimum braids to study graph trees, amphicheirality, unknotting numbers, and periodic tables.

## 1. Introduction

According to Alexander's Theorem [1], any oriented link can be represented as a closed braid. This relationship is not unique in that there is always an infinite number of possible braids for any oriented knot or link. However, by carefully specifying a set of restrictions, it is possible to define a unique braid for any given link.

In this paper we present a new invariant that is called the minimum braid. Since all links can be represented as a braid and since the minimum braid is always unique, the minimum braid is a complete invariant. This invariant is different from a *minimal* braid. A minimal braid has the smallest braid index, i.e., the fewest number of strings or strands. The minimal braid is not unique since many braids for a given link can have the same number of strands. The primary criteria for a minimum braid is that it has the fewest number of braid crossings.

This paper is divided into two main sections. In the first section, a formal definition of minimum braid is presented and the methods used to generated these braids are described. At the end of this section, tables are presented with the minimum braids for all oriented links up to nine crossings and all knots up to ten crossings.

In the second section, a number of applications of minimum braids are described briefly. These applications include trees in graph theory, amphicheiral links, unknotting numbers, and periodic tables.

## 2. Definitions

**Definition 1.** There are four restrictions that are used to define a minimum braid. Among the set of braids for any oriented knot or link, the minimum braid is the one that has the following properties:

- minimum number of braid crossings
- minimum number of braid strands
- minimum braid *universe*
- minimum binary code for *alternating* braid crossings.

These four criteria are listed in descending order of importance when evaluating which element in a set of braids is the minimum braid.

**Definition 2.** A braid *universe* [2] is defined to be an ordered sequence of integers, where the element i represents an unsigned crossing of the i and i+1 braid strands. The braid universe can be considered the projection or "shadow" of a braid word. Any pair of braid universes with the same number of crossings and strands can be sorted by comparing the integers in each sequence. At the first integer which is different, the smaller braid universe is the one with the smaller integer. The minimum braid universe is the first universe generated by this sorting criterion.

A braid universe becomes a braid word when each crossing is assigned a specific overcrossing or undercrossing. A knot or link is formed by joining the ends of each strand of a braid word.

**Definition 3.** A signed crossing or generator in a braid word is defined to be positive when strand i crosses over strand i+1, when the braid is read from top to bottom. This convention is used throughout this paper and regardless of any rotation. Therefore, each of the following crossings are considered to be positive:

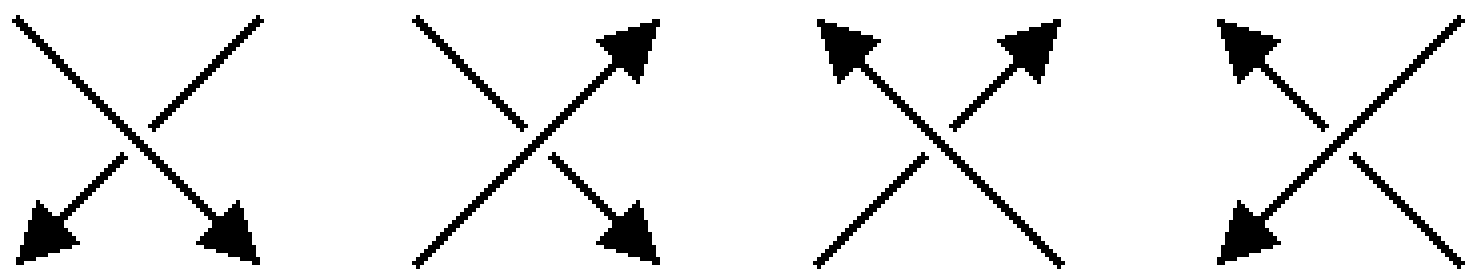

This definition of positive insures that the minimum braids for the majority of knots with ten or fewer crossings will be consistent with the orientations in Rolfsen's *Knots and Links*. [3] It is opposite the convention normally used to define wirthe and polynomials that are based on the skein relationship.

In a braid word, a positive (negative) crossing of the i and i+1 strands is represented by the ith capital (lower case) letter. Using this convention, the trefoil is left-handed, has the minimum braid AAA and a writhe of positive three..

**Definition 4.**. A braid is defined to be *alternating* [4] if even numbered generators have the opposite sign of odd number generators. With our definition of positive crossings, an alternating braid has all odd numbered braid strands crossing over even numbered braid strands. Therefore, crossings in alternating braids have capitals for the odd letters and lower cases for the even letters.

Even for a non-alternating braid, we shall take the convention that a crossing is said to be *alternating* if it is capital for an odd letter and lower case for an even letter. With our positivity convention, an alternating braid has all alternating crossings. When the number of braid crossings is equal to the number of link crossings, the link will be alternating when the braid is alternating.

A binary code for any braid crossing can be generated by assigning a zero for an alternating crossing and a one for a non-alternating crossing. For any braid universe, each knot or link has a unique binary code. When the sequences of zeros and ones are combined to form a number, the minimum binary code is the smallest number.

With the four criteria that define a minimum braid, there is always a unique minimum for any set of braids. The first three criteria determine the minimum braid universe and the binary code determines the minimum braid. For any knot or orientated links, there is a set of Markov equivalent braids according to Markov's theorem [..]. For this set of braids there is always a unique minimum braid; therefore, the minimum braid is a complete invariant.

**Definition 5.** In order to emphasize the difference between minimum and minimal braids, a new class of links is defined. A link is called "MumMal" when its minimum braid crossing is not minimal. Among the knots and links listed at the end of this section, only one knot is MumMal.

## 3. Generating Minimum Braids

When deciding how to generate the minimum braids for knots and oriented links, two approaches were considered. The first approach would start with a nonminimum braid of a link projection and try to find a sequence of Reidemeister-Markov moves that would reduce this braid to a minimum. This "top down" approach was quickly abandoned because it was too difficult to implement and check. When this project was started, there was no tabulation of minimum braids with which to check the effectiveness of attempts to reduce nonminimum braids..

The approach that was followed is a "bottom up" procedure where a set of braids is first generated in a systematic manner so that the first braid corresponding to a particular link is always the minimum braid. Since most braids are not minimum braids, this approach can create large numbers of braids that are rejected because they represent links for which the minimum braid has already been found.

The first step in finding the minimum braids is to generate all possible braid universes for a given number of strands (s) and braid crossings (c). Since the number of possible crossings is equal to the number of braid strands minus one, the total number of braid universes is $(s-1)^c$.

This number can be greatly reduced by applying a series of filters that eliminate universes that can not generate a minimum braid. The first filter is to eliminate all braid universes that do not begin with the integer 1.

The second filter selects only those universes that have two or more crossings for each strand. If there are no crossings that use a braid strand, the number of braid strand is not a minimum. If there is only a single crossing for any braid strand, the braid strand can be eliminated by using Type I Reidemeister moves, regardless of whether the crossing is over or under.

The third filter eliminates any universe that is not a minimum orientation for the sequence of integers. Since the ends of the braid strands will be joined to form a knot or link, there are c possible starting points to use to form an integer sequence for any braid universe. Furthermore, there are four ways to orient any braid universe. The unsigned crossings can be read from top to bottom or bottom to top. The braid universe can then be turned over and again read from top to bottom or bottom to top. This means that there are 4c possible orientation for any braid universe. Only those braid universes that are in their minimum orientation are accepted.

The fourth filter tests for minimum orientations after trying different switching between adjacent crossings. Switching is allowed whenever the absolute difference between two integers in a braid universe is greater that one. This switching is allowed regardless of what pattern will be assigned for the over and under crossings. After each switch the new braid universe is tested to see if its minimum orientation is less than the initial universe. If switching leads to a new minimum, then the braid universe is rejected.

These four filters can reduced the number of braid universes that must be used to generated braid crossings with specific over- and undercrossing assignments. In the case of five braid strands and ten crossings, there are over one million possible universes. These four filters reduce the number to less than two hundred.

After reducing the number of braid universes to a more manageable size, the universes are then separated into different groups according to the number of components in the links that can then be generated. If the sum (s+c) of the number of braid strands and the number of braid crossings is odd (even), then the number of components is odd (even). Each group is then further divided into subgroups depending if the braid universes generate composite or noncomposite links. For example, when there are five braid strands and ten braid crossings, there are only 30 universes that generate knots, and five of these generate composite knots.

After separating the universes by components and composites, knots and links can be generated by specifying over or under crossing sequences. With c braid crossings, there are $2^c$ sequences that for each universe. Some of these sequences can be dropped since they can not be a minimum braid. The first crossing is assumed to be an overcrossing of the first braid strand over the second. This crossing is labeled with a capital **A**. Reverse or mirror image knots and braids would have an initial lower case **a**.

If the universe has a series of two of more crossings that are the same digit, all crossing in the series are assigned the same over- or undercrossing. A braid sequence with consecutive over- and undercrossing of the same braid strands (i.e. **Aa**) could be simplified by a Type II Reidemeister move.

The number of signed braid sequences can be reduced further by checking for patterns of crossings that are equivalent to other patterns. If any of the equivalent patterns are "less" that the original, then the signed braid sequence can not be a minimum and can be rejected.

After generating braid sequences and removing ones that obviously are not minimums, the next step was to identify each braid presentation with a knot or oriented link. The braids were generated in a way so that the first braid that matched a knot or oriented link would be the minimum braid. This procedure generated braids according to the number of braid crossings, the number of braid strands, the lexicographical ordering of the braid universe, and the length of the initial alternating braid sequence.

For each braid, the Alexander polynomial was calculated and powers of the single variable **t** were factored out. The polynomial was then evaluated by setting **t** equal to ten and the absolute value of the resulting integer was taken. This integer was unique for each knot or links tested whenever the general Alexander polynomial was unique. Different Alexander polynomials for **t** were never found to have the same integer when **t** was set equal to ten.

While the Alexander polynomials are not unique for many knots and links, their integer was useful in filtering out possible braids. For example, the last seven knots found had five braid strands and 14 braid crossings. By knowing what Alexander polynomial integers were needed for the unmatched knots, the number of braids that had to be checked was greatly reduced from the millions of possible braids.

The second step for identifying knots or oriented links was to calculate the HOMFLYPT polynomials. The polynomials for knots were then checked with the tables in Aneziris [6]. The polynomials for the oriented links were checked with the tables in Doll and Hoste [7]. In all cases, these tables appeared to be accurate.

The small numbers of pairs of knots (5) and of oriented links (9) with identical HOMFLYPT polynomials were evaluated manually. The following three tables lists the minimum braids that were found.

Table 1 lists the minimum braids for knots with ten or fewer knot crossings. These knots are sorted by the number of knot crossings, with alternating knots listed before nonalternating knots. The first column (Tag #) is a unique name for each knot that has the number of knot crossings as the first digit(s), followed by a colon and 1 to represent the number of components. For alternating knots the next character is a hyphen (-). For nonalternating knots the next character in **n**. The digits following the character is the number for the knot based on the ordering of their minimum braids.

The second column list the number of braid strands in the first position. The second character is a hyphen (-) for an alternating braid and **n** for a nonalternating braid. The last two digits are the number of braid crossings.

The third column lists the minimum braids. The fourth column is the Alexander polynomial integer. The fifth column is the exponent for the variable **t** that was factored out in the braid Alexander polynomial.

The sixth column is the digital of the Alexander polynomial evaluated at ten. The digital of a number is found by summing the digits of an original number. If this sum is greater than nine, the process is repeated until the sum is less than ten. [8] Notice that all of the digitals for knots are either one or eight. This property will be related to the sum of the coefficients of the Alexander polynomial in the section on periodic tables.

The seventh column is the minimum braid unknotting number. This variable is described in the following section on unknotting numbers.

The final column lists the Rolfsen numbers for the knots. An **R** in the last position indicates that the knot generated by the minimum braid is equal to the inverse of the knot drawn in Rolfsen [3]. The second Perko pair has been eliminated and the ten crossing knots starting with 10n162 have been renumbered, accordingly.

Table 1. Minimum Braids for Knots

| Tag # | S-Cr | Minimum Braid | AP(10) | z | d | u | Rolf # | R |
|---|---|---|---|---|---|---|---|---|
| 3:1-01 | 2-03 | AAA | 91 | 0 | 1 | 1 | 3-01 | |
| 4:1-01 | 3-04 | AbAb | 71 | 0 | 8 | 1 | 4-01 | |
| 5:1-01 | 2-05 | AAAAA | 9091 | 0 | 1 | 2 | 5-01 | |
| 5:1-02 | 3n06 | AAABaB | 172 | 1 | 1 | 1 | 5-02 | |
| 6:1-01 | 3-06 | AAAbAb | 7271 | 0 | 8 | 1 | 6-02 | |
| 6:1-02 | 3-06 | AAbAbb | 7471 | 0 | 1 | 1 | 6-03 | |
| 6:1-03 | 4n07 | AABacBc | 152 | 1 | 8 | 1 | 6-01 | |
| 7:1-01 | 2-07 | AAAAAAA | 909091 | 0 | 1 | 3 | 7-01 | |
| 7:1-02 | 4-07 | AAbACbC | 5651 | 0 | 8 | 1 | 7-06 | |
| 7:1-03 | 4-07 | AbAbCbC | 5851 | 0 | 1 | 1 | 7-07 | R |
| 7:1-04 | 3n08 | AAAAABaB | 17272 | 1 | 1 | 2 | 7-03 | R |
| 7:1-05 | 3n08 | AAAABaBB | 16462 | 1 | 1 | 2 | 7-05 | |
| 7:1-06 | 4n09 | AAABaBCbC | 253 | 2 | 1 | 1 | 7-02 | |
| 7:1-07 | 4n09 | AABaBBCbC | 334 | 2 | 1 | 2 | 7-04 | R |
| 8:1-01 | 3-08 | AAAAAbAb | 727271 | 0 | 8 | 2 | 8-02 | |
| 8:1-02 | 3-08 | AAAAbAbb | 745471 | 0 | 1 | 1 | 8-07 | R |
| 8:1-03 | 3-08 | AAAbAAAb | 735371 | 0 | 8 | 2 | 8-05 | R |
| 8:1-04 | 3-08 | AAAbAAbb | 753571 | 0 | 1 | 2 | 8-10 | R |
| 8:1-05 | 3-08 | AAAbAbbb | 743471 | 0 | 8 | 1 | 8-09 | |
| 8:1-06 | 3-08 | AAbAAbAb | 671761 | 0 | 1 | 2 | 8-16 | |
| 8:1-07 | 3-08 | AAbAbAbb | 669761 | 0 | 8 | 1 | 8-17 | |
| 8:1-08 | 3-08 | AbAbAbAb | 587951 | 0 | 8 | 2 | 8-18 | |
| 8:1-09 | 5-08 | AbACbdCd | 4231 | 0 | 1 | 2 | 8-12 | |
| 8:1-10 | 4n09 | AAAABacBc | 14642 | 1 | 8 | 2 | 8-06 | |
| 8:1-11 | 4n09 | AAAbAbcBc | 15452 | 1 | 8 | 2 | 8-04 | |
| 8:1-12 | 4n09 | AAABaBcBc | 13022 | 1 | 8 | 1 | 8-14 | |
| 8:1-13 | 4n09 | AAABacBcc | 14842 | 1 | 1 | 2 | 8-08 | R |
| 8:1-14 | 4n09 | AAbAbbcBc | 14032 | 1 | 1 | 1 | 8-13 | |
| 8:1-15 | 4n09 | AABaBBcBc | 13832 | 1 | 8 | 1 | 8-11 | |
| 8:1-16 | 4n09 | AAbACBBBC | 23023 | 1 | 1 | 2 | 8-15 | |
| 8:1-17 | 5n10 | AABaBCbdCd | 233 | 2 | 8 | 1 | 8-01 | |
| 8:1-18 | 5n10 | AABacBcdCd | 314 | 2 | 8 | 2 | 8-03 | |
| 8:1n19 | 3n08 | AAAbaaab | 8281 | 1 | 1 | 1 | 8n20 | R |
| 8:1n20 | 3n08 | AAABAAAB | 900991 | 0 | 1 | 3 | 8n19 | R |
| 8:1n21 | 3n08 | AAABaaBB | 6461 | 1 | 8 | 1 | 8n21 | |
| 9:1-01 | 2-09 | AAAAAAAAA | 90909091 | 0 | 1 | 4 | 9-01 | |
| 9:1-02 | 4-09 | AAAAbACbC | 563651 | 0 | 8 | 2 | 9-11 | R |
| 9:1-03 | 4-09 | AAAbAACbC | 571751 | 0 | 8 | 2 | 9-36 | R |
| 9:1-04 | 4-09 | AAAbAbCbC | 598051 | 0 | 1 | 1 | 9-26 | R |
| 9:1-05 | 4-09 | AAAbACbCC | 579851 | 0 | 8 | 2 | 9-20 | |
| 9:1-06 | 4-09 | AAbAbACbC | 524341 | 0 | 1 | 2 | 9-32 | R |
| 9:1-07 | 4-09 | AAbAbbCbC | 596051 | 0 | 8 | 1 | 9-27 | |
| 9:1-08 | 4-09 | AAbAbCbCC | 614251 | 0 | 1 | 2 | 9-31 | |
| 9:1-09 | 4-09 | AAbACbbbC | 587951 | 0 | 8 | 1 | 9-24 | |
| 9:1-10 | 4-09 | AAbACbbCC | 606151 | 0 | 1 | 1 | 9-28 | |
| 9:1-11 | 4-09 | AAbbAbCbC | 604151 | 0 | 8 | 1 | 9-30 | |
| 9:1-12 | 4-09 | AbAbbACbC | 522341 | 0 | 8 | 1 | 9-33 | |
| 9:1-13 | 4-09 | AbAbbbCbC | 581851 | 0 | 1 | 2 | 9-17 | R |
| 9:1-14 | 4-09 | AbAbCbAbC | 538541 | 0 | 8 | 1 | 9-34 | |
| 9:1-15 | 4-09 | AbAbCbbbC | 589951 | 0 | 1 | 1 | 9-22 | |
| 9:1-16 | 4-09 | AbACbACbC | 458731 | 0 | 1 | 2 | 9-40 | |

Table 1. Cont'd.

| | | | | | | | | |
|---|---|---|---|---|---|---|---|---|
| 9:1-17 | 4-09 | AbbCbAbCb | 606151 | 0 | 1 | 2 | 9-29 | R |
| 9:1-18 | 3n10 | AAAAAAABaB | 1727272 | 1 | 1 | 3 | 9-03 | R |
| 9:1-19 | 3n10 | AAAAAABaBB | 1645462 | 1 | 1 | 3 | 9-06 | |
| 9:1-20 | 3n10 | AAAAABaBBB | 1653562 | 1 | 1 | 3 | 9-09 | |
| 9:1-21 | 3n10 | AAAABBaBBB | 1571752 | 1 | 1 | 3 | 9-16 | R |
| 9:1-22 | 5n10 | AAABacBDcD | 11402 | 1 | 8 | 2 | 9-15 | R |
| 9:1-23 | 5n10 | AAbAbcBDcD | 13022 | 1 | 8 | 2 | 9-08 | |
| 9:1-24 | 5n10 | AABaBcBDcD | 10592 | 1 | 8 | 1 | 9-21 | R |
| 9:1-25 | 5n10 | AAbACBBdCd | 19583 | 1 | 8 | 2 | 9-25 | |
| 9:1-26 | 5n10 | AAbACbCDcD | 12212 | 1 | 8 | 1 | 9-12 | |
| 9:1-27 | 5n10 | AABacBcDcD | 12412 | 1 | 1 | 1 | 9-14 | R |
| 9:1-28 | 5n10 | AbAbbcBDcD | 11602 | 1 | 1 | 1 | 9-19 | R |
| 9:1-29 | 4n11 | AAAAABaBCbC | 25453 | 2 | 1 | 2 | 9-04 | |
| 9:1-30 | 4n11 | AAAABaBBCbC | 32014 | 2 | 1 | 3 | 9-13 | R |
| 9:1-31 | 4n11 | AAAABaBCbCC | 23833 | 2 | 1 | 2 | 9-07 | |
| 9:1-32 | 4n11 | AAABaBBBCbC | 31204 | 2 | 1 | 2 | 9-18 | |
| 9:1-33 | 4n11 | AAABaBBCbCC | 30394 | 2 | 1 | 2 | 9-23 | |
| 9:1-34 | 4n11 | AABaBBBBCbC | 32824 | 2 | 1 | 3 | 9-10 | R |
| 9:1-35 | 4n11 | AABBcBaBCCB | 37765 | 2 | 1 | 3 | 9-38 | |
| 9:1-36 | 5n12 | AAABaBCbCDcD | 334 | 3 | 1 | 1 | 9-02 | |
| 9:1-37 | 5n12 | AABaBBCbCDcD | 496 | 3 | 1 | 2 | 9-05 | R |
| 9:1-38 | 5n12 | AAbACbadCbCd | 10792 | 2 | 1 | 2 | 9-37 | |
| 9:1-39 | 5n12 | AABacbADCbCD | 17963 | 2 | 8 | 1 | 9-39 | R |
| 9:1-40 | 5n12 | AABacbbDCbCD | 19783 | 2 | 1 | 2 | 9-41 | |
| 9:1-41 | 5n14 | AABaBBCbbDcBDC | 577 | 4 | 1 | 3 | 9-35 | |
| 9:1n42 | 4n09 | AAAbaaCbC | 8081 | 1 | 8 | 1 | 9n42 | R |
| 9:1n43 | 4n09 | AAABAAcBc | 719171 | 0 | 8 | 2 | 9n43 | R |
| 9:1n44 | 4n09 | AAABaacBc | 6661 | 1 | 1 | 1 | 9n44 | |
| 9:1n45 | 4n09 | AABaBACbC | 4841 | 1 | 8 | 1 | 9n45 | |
| 9:1n46 | 4n09 | AbAbCBaBC | 152 | 2 | 8 | 2 | 9n46 | |
| 9:1n47 | 4n09 | AbAbcbAbc | 655561 | 0 | 1 | 2 | 9n47 | |
| 9:1n48 | 4n11 | AABAAcBaBCC | 24643 | 2 | 1 | 3 | 9n49 | R |
| 9:1n49 | 4n11 | AABaBAcBaBc | 4031 | 2 | 8 | 2 | 9n48 | R |
| 10:1-001 | 3-10 | AAAAAAAbAb | 72727271 | 0 | 8 | 3 | 10-002 | |
| 10:1-002 | 3-10 | AAAAAAbAbb | 74545471 | 0 | 1 | 2 | 10-005 | R |
| 10:1-003 | 3-10 | AAAAAbAAAb | 73545371 | 0 | 8 | 3 | 10-046 | R |
| 10:1-004 | 3-10 | AAAAAbAAbb | 75363571 | 0 | 1 | 3 | 10-047 | R |
| 10:1-005 | 3-10 | AAAAAbAbbb | 74363471 | 0 | 8 | 1 | 10-009 | R |
| 10:1-006 | 3-10 | AAAAbAAAbb | 75282571 | 0 | 1 | 2 | 10-062 | R |
| 10:1-007 | 3-10 | AAAAbAAbAb | 67100761 | 0 | 1 | 2 | 10-085 | |
| 10:1-008 | 3-10 | AAAAbAbAbb | 66918761 | 0 | 8 | 1 | 10-082 | |
| 10:1-009 | 3-10 | AAAAbAbbbb | 74383471 | 0 | 1 | 1 | 10-017 | |
| 10:1-010 | 3-10 | AAAAbbAbbb | 75201571 | 0 | 1 | 2 | 10-048 | |
| 10:1-011 | 3-10 | AAAbAAAbbb | 75100571 | 0 | 8 | 2 | 10-064 | R |
| 10:1-012 | 3-10 | AAAbAAbAAb | 67918861 | 0 | 1 | 3 | 10-100 | |
| 10:1-013 | 3-10 | AAAbAAbbAb | 67736861 | 0 | 8 | 2 | 10-094 | R |
| 10:1-014 | 3-10 | AAAbAbAAbb | 67655861 | 0 | 8 | 2 | 10-106 | R |
| 10:1-015 | 3-10 | AAAbAbAbAb | 59474051 | 0 | 8 | 2 | 10-112 | |
| 10:1-016 | 3-10 | AAAbAbbAbb | 67756861 | 0 | 1 | 1 | 10-091 | |
| 10:1-017 | 3-10 | AAAbbAAbbb | 75938671 | 0 | 1 | 3 | 10-079 | |
| 10:1-018 | 3-10 | AAAbbAbAbb | 67675861 | 0 | 1 | 1 | 10-104 | |
| 10:1-019 | 3-10 | AAbAAbAbAb | 60292151 | 0 | 8 | 2 | 10-116 | |
| 10:1-020 | 3-10 | AAbAAbbAbb | 68574961 | 0 | 1 | 2 | 10-099 | |

Table 1. Cont'd.

| | | | | | | | | |
|---|---|---|---|---|---|---|---|---|
| 10:1-021 | 3-10 | AAbAbAbbAb | 60312151 | 0 | 1 | 1 | 10-118 | R |
| 10:1-022 | 3-10 | AAbAbbAAbb | 68493961 | 0 | 1 | 2 | 10-109 | |
| 10:1-023 | 3-10 | AbAbAbAbAb | 52867441 | 0 | 1 | 2 | 10-123 | |
| 10:1-024 | 5-10 | AAAbACbdCd | 434431 | 0 | 1 | 2 | 10-029 | |
| 10:1-025 | 5-10 | AAbAbCbdCd | 464831 | 0 | 8 | 1 | 10-042 | |
| 10:1-026 | 5-10 | AAbACbbdCd | 456731 | 0 | 8 | 1 | 10-071 | |
| 10:1-027 | 5-10 | AAbACbCdCd | 466831 | 0 | 1 | 1 | 10-044 | |
| 10:1-028 | 5-10 | AAbACbdCdd | 448631 | 0 | 8 | 2 | 10-043 | |
| 10:1-029 | 5-10 | AbAbbCbdCd | 450631 | 0 | 1 | 2 | 10-041 | R |
| 10:1-030 | 5-10 | AbAbCbbdCd | 458731 | 0 | 1 | 1 | 10-059 | |
| 10:1-031 | 5-10 | AbAbCbCdCd | 481031 | 0 | 8 | 2 | 10-045 | |
| 10:1-032 | 5-10 | AbACbbbdCd | 442531 | 0 | 1 | 2 | 10-070 | |
| 10:1-033 | 5-10 | AbACbCbdCd | 407321 | 0 | 8 | 1 | 10-088 | |
| 10:1-034 | 4n11 | AAAAAABacBc | 1463642 | 1 | 8 | 3 | 10-006 | |
| 10:1-035 | 4n11 | AAAAAbAbcBc | 1545452 | 1 | 8 | 2 | 10-008 | |
| 10:1-036 | 4n11 | AAAAABaBcBc | 1308122 | 1 | 8 | 2 | 10-014 | R |
| 10:1-037 | 4n11 | AAAAABacBcc | 1489942 | 1 | 1 | 2 | 10-012 | R |
| 10:1-038 | 4n11 | AAAAbAbbcBc | 1400032 | 1 | 1 | 2 | 10-019 | |
| 10:1-039 | 4n11 | AAAABaBBcBc | 1324322 | 1 | 8 | 2 | 10-025 | |
| 10:1-040 | 4n11 | AAAAbAbcBcc | 1481842 | 1 | 1 | 2 | 10-015 | R |
| 10:1-041 | 4n11 | AAAABaBcBcc | 1342522 | 1 | 1 | 1 | 10-027 | |
| 10:1-042 | 4n11 | AAAAbACBBBC | 2308123 | 1 | 1 | 3 | 10-049 | |
| 10:1-043 | 4n11 | AAAABacBBBc | 1406132 | 1 | 8 | 3 | 10-076 | R |
| 10:1-044 | 4n11 | AAAABacBBcc | 1424332 | 1 | 1 | 3 | 10-077 | R |
| 10:1-045 | 4n11 | AAAABacBccc | 1487942 | 1 | 8 | 2 | 10-022 | R |
| 10:1-046 | 4n11 | AAAABBaBcBc | 1242512 | 1 | 8 | 2 | 10-072 | R |
| 10:1-047 | 4n11 | AAAbAAAbcBc | 1553552 | 1 | 8 | 3 | 10-061 | R |
| 10:1-048 | 4n11 | AAAbAAbbcBc | 1416232 | 1 | 1 | 2 | 10-052 | R |
| 10:1-049 | 4n11 | AAAbAAbcBcc | 1489942 | 1 | 1 | 3 | 10-054 | R |
| 10:1-050 | 4n11 | AAAbAACBBBC | 2234413 | 1 | 1 | 3 | 10-080 | |
| 10:1-051 | 4n11 | AAAbAbbbcBc | 1414232 | 1 | 8 | 1 | 10-026 | |
| 10:1-052 | 4n11 | AAABaBBBcBc | 1316222 | 1 | 8 | 2 | 10-039 | |
| 10:1-053 | 4n11 | AAAbAbbcBcc | 1332422 | 1 | 8 | 1 | 10-032 | R |
| 10:1-054 | 4n11 | AAABaBBcBcc | 1350622 | 1 | 1 | 2 | 10-040 | R |
| 10:1-055 | 4n11 | AAABaBcBBBc | 1324322 | 1 | 8 | 2 | 10-056 | R |
| 10:1-056 | 4n11 | AAABaBcBBcc | 1358722 | 1 | 1 | 2 | 10-057 | R |
| 10:1-057 | 4n11 | AAAbACBBBCC | 2242513 | 1 | 1 | 3 | 10-066 | |
| 10:1-058 | 4n11 | AAABacBBcBc | 1276912 | 1 | 1 | 1 | 10-084 | R |
| 10:1-059 | 4n11 | AAABacBcBcc | 1258712 | 1 | 8 | 2 | 10-087 | R |
| 10:1-060 | 4n11 | AAABBcBaBcB | 1176902 | 1 | 8 | 2 | 10-092 | R |
| 10:1-061 | 4n11 | AAABcBaBcBc | 1129492 | 1 | 1 | 1 | 10-113 | R |
| 10:1-062 | 4n11 | AAbAAbAbcBc | 1334422 | 1 | 1 | 2 | 10-093 | |
| 10:1-063 | 4n11 | AAbAACbAbcc | 1326322 | 1 | 1 | 2 | 10-108 | R |
| 10:1-064 | 4n11 | AAbAbAbbcBc | 1266812 | 1 | 8 | 2 | 10-086 | |
| 10:1-065 | 4n11 | AAbAbbbbcBc | 1416232 | 1 | 1 | 1 | 10-023 | |
| 10:1-066 | 4n11 | AABaBBBBcBc | 1381832 | 1 | 8 | 2 | 10-021 | |
| 10:1-067 | 4n11 | AABaBBcBBBc | 1398032 | 1 | 8 | 2 | 10-050 | R |
| 10:1-068 | 4n11 | AABaBBcBBcc | 1432432 | 1 | 1 | 3 | 10-051 | R |
| 10:1-069 | 4n11 | AAbAbcBAcbb | 1348622 | 1 | 8 | 2 | 10-090 | |
| 10:1-070 | 4n11 | AABaBcBBBcc | 1424332 | 1 | 1 | 2 | 10-065 | R |
| 10:1-071 | 4n11 | AABaBcBBcBc | 1268812 | 1 | 1 | 2 | 10-083 | R |
| 10:1-072 | 4n11 | AAbACbAbbcc | 1340522 | 1 | 8 | 1 | 10-102 | |
| 10:1-073 | 4n11 | AAbACbAbccb | 1201202 | 1 | 8 | 1 | 10-119 | |

Table 1. Cont'd.

| | | | | | | | | |
|---|---|---|---|---|---|---|---|---|
| 10:1-074 | 4n11 | AABacBBcBBc | 1350622 | 1 | 1 | 3 | 10-103 | |
| 10:1-075 | 4n11 | AABacBcBcBc | 1185002 | 1 | 8 | 1 | 10-114 | |
| 10:1-076 | 4n11 | AABBcBaBBcB | 1258712 | 1 | 8 | 2 | 10-098 | |
| 10:1-077 | 4n11 | AABBcBaBcBc | 1211302 | 1 | 1 | 2 | 10-117 | R |
| 10:1-078 | 4n11 | AAbbCbAbccb | 1285012 | 1 | 1 | 1 | 10-095 | |
| 10:1-079 | 4n11 | AABBcBBaBcB | 1250612 | 1 | 8 | 2 | 10-111 | R |
| 10:1-080 | 4n11 | AABcBaBcBcB | 1137592 | 1 | 1 | 2 | 10-121 | |
| 10:1-081 | 4n11 | AABcBacBcBc | 1111292 | 1 | 8 | 2 | 10-122 | R |
| 10:1-082 | 6n11 | AABacBDceDe | 9172 | 1 | 1 | 2 | 10-013 | |
| 10:1-083 | 6n11 | AbAbcBDceDe | 9982 | 1 | 1 | 2 | 10-035 | |
| 10:1-084 | 6n11 | AbACbdccEdE | 16543 | 1 | 1 | 2 | 10-058 | R |
| 10:1-085 | 5n12 | AAAABaBCbdCd | 22013 | 2 | 8 | 2 | 10-020 | |
| 10:1-086 | 5n12 | AAAABacBcdCd | 30194 | 2 | 8 | 3 | 10-011 | |
| 10:1-087 | 5n12 | AAABaBBCbdCd | 26954 | 2 | 8 | 2 | 10-038 | |
| 10:1-088 | 5n12 | AAAbAbcBcdCd | 23633 | 2 | 8 | 2 | 10-004 | |
| 10:1-089 | 5n12 | AAABaBCbCdCd | 18773 | 2 | 8 | 2 | 10-036 | |
| 10:1-090 | 5n12 | AAABaBcBcdCd | 27764 | 2 | 8 | 1 | 10-018 | |
| 10:1-091 | 5n12 | AAABaBCbdCdd | 22213 | 2 | 1 | 2 | 10-034 | R |
| 10:1-092 | 5n12 | AAABacBccdCd | 27964 | 2 | 1 | 1 | 10-031 | |
| 10:1-093 | 5n12 | AAABacBcdCdd | 28774 | 2 | 1 | 2 | 10-037 | |
| 10:1-094 | 5n12 | AAABacBDCCCD | 36955 | 2 | 1 | 2 | 10-055 | |
| 10:1-095 | 5n12 | AABaBAcBcDcD | 474931 | 1 | 1 | 1 | 10-073 | |
| 10:1-096 | 5n12 | AABaBAcBDcDD | 440531 | 1 | 8 | 2 | 10-078 | |
| 10:1-097 | 5n12 | AABaBBBCbdCd | 27764 | 2 | 8 | 2 | 10-024 | |
| 10:1-098 | 5n12 | AAbAbbcBcdCd | 20593 | 2 | 1 | 1 | 10-010 | |
| 10:1-099 | 5n12 | AABaBBCbCdCd | 25334 | 2 | 8 | 1 | 10-030 | |
| 10:1-100 | 5n12 | AABaBBcBcdCd | 29384 | 2 | 8 | 2 | 10-016 | R |
| 10:1-101 | 5n12 | AABaBBCbdCdd | 28774 | 2 | 1 | 2 | 10-028 | R |
| 10:1-102 | 5n12 | AABaBCbCCdCd | 20393 | 2 | 8 | 1 | 10-007 | |
| 10:1-103 | 5n12 | AABaBcBccdCd | 26344 | 2 | 1 | 1 | 10-033 | |
| 10:1-104 | 5n12 | AABaBcBDCCCD | 44326 | 2 | 1 | 3 | 10-053 | |
| 10:1-105 | 5n12 | AABacBADcBcD | 483031 | 1 | 1 | 2 | 10-069 | R |
| 10:1-106 | 5n12 | AAbACBBBCDcD | 37765 | 2 | 1 | 2 | 10-063 | |
| 10:1-107 | 5n12 | AAbACBBdcBcd | 393121 | 1 | 1 | 2 | 10-105 | R |
| 10:1-108 | 5n12 | AAbAcbbDcBcD | 391121 | 1 | 8 | 1 | 10-107 | |
| 10:1-109 | 5n12 | AAbACBBdcccd | 374921 | 1 | 8 | 2 | 10-081 | R |
| 10:1-110 | 5n12 | AbAbbCbCBDcD | 472931 | 1 | 8 | 1 | 10-060 | |
| 10:1-111 | 5n12 | AbAbcBADCbCD | 409321 | 1 | 1 | 2 | 10-089 | |
| 10:1-112 | 5n12 | AbAbCbbDcBDC | 464831 | 1 | 8 | 2 | 10-075 | R |
| 10:1-113 | 5n12 | AbaCbaCdCbCd | 489131 | 1 | 8 | 2 | 10-096 | |
| 10:1-114 | 5n12 | AbACBBBdcBcd | 376921 | 1 | 1 | 2 | 10-110 | |
| 10:1-115 | 5n12 | AbACBBdcBccd | 325511 | 1 | 8 | 2 | 10-115 | R |
| 10:1-116 | 6n13 | AABaBCbCDceDe | 314 | 3 | 8 | 1 | 10-001 | |
| 10:1-117 | 6n13 | AABaBCbdCdeDe | 476 | 3 | 8 | 2 | 10-003 | |
| 10:1-118 | 5n14 | AAABaBCbbdCBdC | 26144 | 3 | 8 | 2 | 10-067 | |
| 10:1-119 | 5n14 | AAABaCbACBBDcD | 51697 | 3 | 1 | 3 | 10-101 | R |
| 10:1-120 | 5n14 | AABaBAcBaBCdCd | 31085 | 3 | 8 | 2 | 10-097 | R |
| 10:1-121 | 5n14 | AAbAbbcBBdCbdc | 27964 | 3 | 1 | 2 | 10-068 | R |
| 10:1-122 | 5n14 | AABaBBCbbdCBdC | 26144 | 3 | 8 | 2 | 10-074 | |
| 10:1-123 | 5n14 | AABacbADCBBCCD | 57448 | 3 | 1 | 3 | 10-120 | |
| 10:1n124 | 3n10 | AAAAAbaaab | 819181 | 1 | 1 | 2 | 10n125 | R |
| 10:1n125 | 3n10 | AAAAABAAAB | 90090991 | 0 | 1 | 4 | 10n124 | R |
| 10:1n126 | 3n10 | AAAAABaaaB | 835381 | 1 | 1 | 2 | 10n126 | |

Table 1. Cont'd.

| | | | | | | | | |
|---|---|---|---|---|---|---|---|---|
| 10:1n127 | 3n10 | AAAAABaaBB | 653561 | 1 | 8 | 2 | 10n127 | |
| 10:1n128 | 3n10 | AAAAbaaabb | 735371 | 1 | 8 | 1 | 10n141 | R |
| 10:1n129 | 3n10 | AAAABAAABB | 90171991 | 0 | 1 | 4 | 10n139 | R |
| 10:1n130 | 3n10 | AAAABaaaBB | 753571 | 1 | 1 | 1 | 10n143 | |
| 10:1n131 | 3n10 | AAAABaaBaB | 761671 | 1 | 1 | 2 | 10n148 | |
| 10:1n132 | 3n10 | AAAABaBaBB | 579851 | 1 | 8 | 2 | 10n149 | |
| 10:1n133 | 3n10 | AAABaaBaaB | 743471 | 1 | 8 | 2 | 10n155 | R |
| 10:1n134 | 3n10 | AAABaBAABB | 982801 | 1 | 1 | 3 | 10n161 | |
| 10:1n135 | 3n10 | AAABaBaaBB | 679861 | 1 | 1 | 1 | 10n159 | |
| 10:1n136 | 3n10 | AAABBAABBB | 89353891 | 0 | 1 | 4 | 10n152 | |
| 10:1n137 | 3n10 | AAABBaBaBB | 498041 | 1 | 8 | 2 | 10n157 | R |
| 10:1n138 | 5n10 | AbAbCBBdCd | 5041 | 1 | 1 | 1 | 10n137 | |
| 10:1n139 | 5n10 | AbAbcbbDcD | 573751 | 0 | 1 | 2 | 10n138 | |
| 10:1n140 | 5n10 | AbAbcBBDcD | 6461 | 1 | 8 | 1 | 10n136 | R |
| 10:1n141 | 4n11 | AAAbaaabcBc | 8281 | 2 | 1 | 2 | 10n140 | R |
| 10:1n142 | 4n11 | AAABAAABCbC | 1719172 | 1 | 1 | 3 | 10n142 | R |
| 10:1n143 | 4n11 | AAAbaabbcBc | 16462 | 2 | 1 | 2 | 10n130 | R |
| 10:1n144 | 4n11 | AAABAABBCbC | 1711072 | 1 | 1 | 3 | 10n128 | R |
| 10:1n145 | 4n11 | AAABaaBBCbC | 13022 | 2 | 8 | 1 | 10n131 | |
| 10:1n146 | 4n11 | AAAbaabcBcc | 9091 | 2 | 1 | 1 | 10n132 | R |
| 10:1n147 | 4n11 | AAABAABCbCC | 1637362 | 1 | 1 | 3 | 10n134 | R |
| 10:1n148 | 4n11 | AAABaaBCbCC | 5651 | 2 | 8 | 1 | 10n133 | |
| 10:1n149 | 4n11 | AAAbaaCBaBC | 671761 | 1 | 1 | 1 | 10n156 | |
| 10:1n150 | 4n11 | AAABAAcBaBc | 637361 | 1 | 8 | 2 | 10n160 | R |
| 10:1n151 | 4n11 | AAABaacbAbc | 685961 | 1 | 8 | 2 | 10n158 | |
| 10:1n152 | 4n11 | AAAbAACbaCB | 653561 | 1 | 8 | 2 | 10n150 | R |
| 10:1n153 | 4n11 | AAAbaaCbaCb | 14842 | 2 | 1 | 1 | 10n129 | |
| 10:1n154 | 4n11 | AAABaaCbACb | 687961 | 1 | 1 | 2 | 10n151 | R |
| 10:1n155 | 4n11 | AAABAAcbbbc | 892891 | 1 | 1 | 2 | 10n153 | |
| 10:1n156 | 4n11 | AAAbAbcBaBc | 13832 | 2 | 8 | 1 | 10n147 | R |
| 10:1n157 | 4n11 | AAABaBcbbbc | 22213 | 2 | 1 | 2 | 10n135 | R |
| 10:1n158 | 4n11 | AABaaBBAcBc | 22013 | 2 | 8 | 2 | 10n162 | |
| 10:1n159 | 4n11 | AAbaaCBaBBC | 606151 | 1 | 1 | 2 | 10n163 | R |
| 10:1n160 | 4n11 | AAbAbaCbAbC | 13222 | 2 | 1 | 1 | 10n146 | |
| 10:1n161 | 4n11 | AABaBACBaBC | 10711 | 2 | 1 | 2 | 10n145 | |
| 10:1n162 | 4n11 | AABaBAcBAcb | 21203 | 2 | 8 | 2 | 10n144 | |
| 10:1n163 | 4n11 | AABaBACBBBC | 966601 | 1 | 1 | 3 | 10n154 | R |
| 10:1n164 | 4n11 | AAbAbbcBaBc | 20593 | 2 | 1 | 1 | 10n164 | R |
| 10:1n165 | 4n11 | AABacBaBCCB | 11402 | 2 | 8 | 2 | 10n165 | R |

Table 2 lists the minimum braids for oriented links with nine or fewer crossings. The columns are similar to Table 1 with a couple differences. The sixth column is the digital of the Alexander polynomial divided by $(t-1)^{k-1}$, where **k** is the number of components and the polynomial is evaluated at t=10. The seventh column is the minimum braid unlinking number or the minimum number of braid crossings that need to be switched to reduce the links to **k** unknots. The first column of Tag #'s includes a lower case suffix to identify different oriented links that have the same nonoriented link. The final column lists the identifiers used by Doll and Hoste [7] for oriented links. One of the links (9:2n59a) is an orientation that was not listed by Doll and Hoste.

Table 3 lists the minimum braids for some of the trivial knots and links. These trivials include composite knots and/or links and links formed by adding simple loops to knots or oriented links. In the first columns, the Tag #'s have a lower case **c** following the number of components for composites and an "**.**" for "loops". When there are loops or composites with different orientations, there is a lower case letter suffix used to differentiate the links. These is also one example (8:2_01) of a link formed by joining two knots with a pair of crossings.

It is interesting to compare these minimum braids with braids that others have published for knots through ten crossings. The two lists that were compared were the braids generated by Jones [9] and by Chalcraft [10]. Table 4 summarizes the differences between the minimum braids and the Jones braids.

Table 4. Differences for Jones' braids

| | | Differences | | | | | | |
|---|---|---|---|---|---|---|---|---|
| | | Strands | | Crossings | | | | |
| Strands | Knots | -1 | 0 | 0 | 1 | 2 | 3 | 4 |
| 2 | 4 | - | 4 | 4 | - | - | - | - |
| 3 | 58 | - | 58 | 58 | - | - | - | - |
| 4 | 116 | - | 116 | 90 | - | 26 | - | - |
| 5 | 66 | 1 | 65 | 30 | 1 | 22 | - | 13 |
| 6 | 5 | - | 5 | 2 | - | 1 | - | 2 |
| Totals | 249 | 1 | 248 | 184 | 1 | 49 | 0 | 15 |

Probably the most interesting number in Table 4 is the one in the first column of the differences in strands. The minimum braid for this knot uses five strands while the Jones braid uses just four. This knot is 10:1n140 ($10_{136}$) and has the ten crossing minimum braid AbAbcBBDcD. The Jones braid has four strands but 11 crossings. This is the first example of a MumMal knot.

Table 2. Minimum Braids for Links

| Tag # | S-Cr | Minimum Braid | AP(10) | z | d | u | D&H # |
|---|---|---|---|---|---|---|---|
| 2:2-01 | 2-02 | AA | 9 | 0 | 1 | 1 | 2-2-01 |
| 4:2-01a | 2-04 | AAAA | 909 | 0 | 2 | 2 | 4-2-01+- |
| 4:2-01b | 3n05 | AABaB | 18 | 1 | 2 | 2 | 4-2-01 |
| 5:2-01 | 3-05 | AAbAb | 729 | 0 | 9 | 1 | 5-2-01 |
| 6:2-01a | 2-06 | AAAAAA | 90909 | 0 | 3 | 3 | 6-2-01 |
| 6:2-01b | 4n08 | AABaBCbC | 27 | 2 | 3 | 3 | 6-2-01+- |
| 6:2-02a | 4-06 | AbACbC | 549 | 0 | 7 | 2 | 6-2-03+- |
| 6:2-02b | 3n07 | AAABaBB | 1638 | 1 | 2 | 2 | 6-2-03 |
| 6:2-03 | 3n07 | AAAABaB | 1728 | 1 | 3 | 3 | 6-2-02 |
| 6:3-01a | 3-06 | AAbAAb | 7371 | 0 | 1 | 3 | 6-3-01++- |
| 6:3-01b | 4n07 | AbACBBC | 243 | 1 | 3 | 3 | 6-3-01 |
| 6:3-02 | 3-06 | AbAbAb | 6561 | 0 | 9 | 2 | 6-3-02 |
| 6:3n03a | 3n06 | AAbaab | 81 | 1 | 1 | 3 | 6-3-03 |
| 6:3n03b | 3n06 | AABAAB | 8991 | 0 | 3 | 3 | 6-3-03++- |
| 7:2-01a | 3-07 | AAAAbAb | 72729 | 0 | 8 | 2 | 7-2-01 |
| 7:2-01b | 4n08 | AAbAbcBc | 1548 | 1 | 1 | 2 | 7-2-01+- |
| 7:2-02 | 3-07 | AAAbAAb | 73629 | 0 | 9 | 2 | 7-2-04 |
| 7:2-03a | 3-07 | AAAbAbb | 74529 | 0 | 1 | 1 | 7-2-02+- |
| 7:2-03b | 4n08 | AABaBcBc | 1368 | 1 | 8 | 1 | 7-2-02 |
| 7:2-04a | 3-07 | AAbAAbb | 75429 | 0 | 2 | 3 | 7-2-05+- |
| 7:2-04b | 4n08 | AAbACBBC | 2367 | 1 | 2 | 3 | 7-2-05 |
| 7:2-05 | 3-07 | AAbAbAb | 66339 | 0 | 9 | 2 | 7-2-06 |
| 7:2-06 | 4n08 | AAABacBc | 1458 | 1 | 9 | 2 | 7-2-03 |
| 7:2n07a | 3n07 | AAAbaab | 909 | 1 | 2 | 3 | 7-2-07 |
| 7:2n07b | 3n07 | AAABAAB | 90009 | 0 | 2 | 3 | 7-2-07+- |
| 7:2n08 | 3n07 | AAABaaB | 729 | 1 | 9 | 1 | 7-2-08 |
| 7:3-01a | 4-07 | AbACbbC | 5751 | 0 | 8 | 3 | 7-3-01 |
| 7:3-01b | 3n08 | AAABBaBB | 15552 | 1 | 3 | 3 | 7-3-01+-- |
| 8:2-01a | 2-08 | AAAAAAAA | 9090909 | 0 | 4 | 4 | 8-2-01 |
| 8:2-01b | 5n11 | AABaBCbCDcD | 36 | 3 | 4 | 4 | 8-2-01+- |
| 8:2-02a | 4-08 | AAAbACbC | 56349 | 0 | 6 | 3 | 8-2-05+- |
| 8:2-02b | 4n10 | AAABaBBCbC | 3186 | 2 | 3 | 3 | 8-2-05 |
| 8:2-03a | 4-08 | AAbAACbC | 57249 | 0 | 7 | 3 | 8-2-09+- |
| 8:2-03b | 5n09 | AbACBBdCd | 2007 | 1 | 7 | 3 | 8-2-09 |
| 8:2-04 | 4-08 | AAbAbCbC | 59949 | 0 | 1 | 2 | 8-2-08 |
| 8:2-05 | 4-08 | AAbACbbC | 59049 | 0 | 9 | 1 | 8-2-12 |
| 8:2-06a | 4-08 | AAbACbCC | 58149 | 0 | 8 | 2 | 8-2-07+- |
| 8:2-06b | 4-08 | AbAbbCbC | 58149 | 0 | 8 | 2 | 8-2-07 |
| 8:2-07 | 4-08 | AbAbACbC | 51759 | 0 | 9 | 2 | 8-2-13 |
| 8:2-08 | 4-08 | AbAbCbbC | 59049 | 0 | 9 | 1 | 8-2-10 |
| 8:2-09a | 4-08 | AbACbbbC | 57249 | 0 | 7 | 3 | 8-2-11+- |
| 8:2-09b | 3n09 | AAAABBaBB | 156348 | 1 | 2 | 3 | 8-2-11 |
| 8:2-10a | 4-08 | AbCbAbCb | 60849 | 0 | 2 | 3 | 8-2-14+- |
| 8:2-10b | 4n10 | AABcBaBCCB | 3825 | 2 | 2 | 3 | 8-2-14 |
| 8:2-11a | 3n09 | AAAAAABaB | 172728 | 1 | 4 | 4 | 8-2-02 |
| 8:2-11b | 4n10 | AAAABaBCbC | 2547 | 2 | 4 | 4 | 8-2-02+- |
| 8:2-12a | 3n09 | AAAAABaBB | 164538 | 1 | 3 | 3 | 8-2-03 |
| 8:2-12b | 5n09 | AABacBDcD | 1188 | 1 | 6 | 3 | 8-2-03+- |
| 8:2-13a | 3n09 | AAAABaBBB | 165438 | 1 | 4 | 4 | 8-2-04+- |
| 8:2-13b | 4n10 | AABaBBBCbC | 3276 | 2 | 4 | 4 | 8-2-04 |
| 8:2-14a | 5n09 | AbAbcBDcD | 1278 | 1 | 7 | 2 | 8-2-06+- |

Table 2. Cont'd.

| | | | | | | | |
|---|---|---|---|---|---|---|---|
| 8:2-14b | 4n10 | AAABaBCbCC | 2367 | 2 | 2 | 2 | 8-2-06 |
| 8:2n15 | 4n08 | AAbaaCbC | 729 | 1 | 9 | 1 | 8-2-15 |
| 8:2n16a | 4n08 | AABAAcBc | 71829 | 0 | 7 | 3 | 8-2-16+- |
| 8:2n16b | 4n08 | ABaBACbC | 549 | 1 | 7 | 3 | 8-2-16 |
| 8:3-01a | 3-08 | AAAAbAAb | 736371 | 0 | 1 | 4 | 8-3-01 |
| 8:3-01b | 4n09 | AAbAAbcBc | 15552 | 1 | 3 | 4 | 8-3-01+-+ |
| 8:3-01c | 5n10 | AABacBDCCD | 405 | 2 | 5 | 4 | 8-3-01++- |
| 8:3-02a | 3-08 | AAAbAbAb | 662661 | 0 | 9 | 3 | 8-3-05 |
| 8:3-02b | 4n09 | AABacBcBc | 13122 | 1 | 9 | 3 | 8-3-05++- |
| 8:3-03a | 3-08 | AAAbbAbb | 752571 | 0 | 3 | 4 | 8-3-02+-- |
| 8:3-03b | 4n09 | AAAbACBBC | 23733 | 1 | 5 | 4 | 8-3-02 |
| 8:3-03c | 4n09 | AABaBcBBc | 13932 | 1 | 1 | 4 | 8-3-02+-+ |
| 8:3-04a | 3-08 | AAbAbbAb | 670761 | 0 | 1 | 2 | 8-3-06 |
| 8:3-04b | 4n09 | AABcBaBcB | 12312 | 1 | 8 | 2 | 8-3-06++- |
| 8:3-05a | 3-08 | AAbbAAbb | 760671 | 0 | 4 | 4 | 8-3-04++- |
| 8:3-05b | 4n09 | AAbAACBBC | 22923 | 1 | 4 | 4 | 8-3-04+-+ |
| 8:3-05c | 5n10 | AbACBBdccd | 3321 | 1 | 5 | 4 | 8-3-04 |
| 8:3-06a | 4n09 | AAABacBBc | 13932 | 1 | 1 | 3 | 8-3-03++- |
| 8:3-06b | 5n10 | ABaBAcBDcD | 4131 | 1 | 6 | 3 | 8-3-03 |
| 8:3n07a | 3n08 | AAAAbaab | 8991 | 1 | 3 | 4 | 8-3-07 |
| 8:3n07b | 3n08 | AAAABAAB | 899991 | 0 | 5 | 4 | 8-3-07+-- |
| 8:3n07c | 4n09 | AAbaabcBc | 81 | 2 | 1 | 4 | 8-3-07++- |
| 8:3n08a | 3n08 | AAAABaaB | 7371 | 1 | 1 | 4 | 8-3-08 |
| 8:3n08b | 4n09 | AABAABCbC | 17172 | 1 | 5 | 4 | 8-3-08++- |
| 8:3n08c | 4n09 | AABaBcbbc | 243 | 2 | 3 | 4 | 8-3-08+-+ |
| 8:3n09a | 3n08 | AAABaBaB | 6561 | 1 | 9 | 2 | 8-3-09 |
| 8:3n09b | 4n09 | AAbaaCBBC | 6561 | 1 | 9 | 2 | 8-3-09++- |
| 8:3n10a | 3n08 | AABBAABB | 891891 | 0 | 4 | 4 | 8-3-10++- |
| 8:3n10b | 4n09 | AABAAcbbc | 9801 | 1 | 4 | 4 | 8-3-10+-+ |
| 8:3n10c | 4n09 | ABaBACBBC | 9801 | 1 | 4 | 4 | 8-3-10 |
| 8:4-01a | 4-08 | AbbACbbC | 59049 | 0 | 9 | 4 | 8-4-01++-- |
| 8:4-01b | 4n10 | AAAbbacBBc | 59049 | 1 | 9 | 4 | 8-4-01+++- |
| 8:4-01c | 4n10 | AAABBaCBBC | 147258 | 1 | 4 | 4 | 8-4-01+--+ |
| 8:4-01d | 5n11 | ABaBAcBDCCD | 2916 | 2 | 4 | 4 | 8-4-01 |
| 8:4n02a | 4n08 | AbbACBBC | 1458 | 1 | 2 | 4 | 8-4-02 |
| 8:4n02b | 4n08 | AbbAcbbc | 73629 | 0 | 2 | 4 | 8-4-02++-- |
| 8:4n03a | 4n08 | AbbAcBBc | 0 | 2 | 9 | 4 | 8-4-03+++- |
| 8:4n03b | 4n08 | ABBACBBC | 88209 | 0 | 4 | 4 | 8-4-03++-- |
| 8:4n03c | 4n10 | AAAbbaCbbC | 0 | 3 | 9 | 4 | 8-4-03 |
| 9:2-01a | 3-09 | AAAAAAbAb | 7272729 | 0 | 7 | 3 | 9-2-01 |
| 9:2-01b | 5n11 | AAbAbcBcdCd | 2367 | 2 | 2 | 3 | 9-2-01+- |
| 9:2-02 | 3-09 | AAAAAbAAb | 7363629 | 0 | 9 | 3 | 9-2-13 |
| 9:2-03a | 3-09 | AAAAAbAbb | 7454529 | 0 | 2 | 2 | 9-2-02 |
| 9:2-03b | 5n11 | AABaBCbCdCd | 2007 | 2 | 7 | 2 | 9-2-02+- |
| 9:2-04a | 3-09 | AAAAbAAAb | 7354629 | 0 | 8 | 3 | 9-2-19 |
| 9:2-04b | 4n10 | AAAbAAbcBc | 155448 | 1 | 1 | 3 | 9-2-19+- |
| 9:2-05a | 3-09 | AAAAbAAbb | 7536429 | 0 | 3 | 4 | 9-2-20 |
| 9:2-05b | 5n11 | AAbACBBCDcD | 3915 | 2 | 3 | 4 | 9-2-20+- |
| 9:2-06 | 3-09 | AAAAbAbAb | 6627339 | 0 | 9 | 3 | 9-2-31 |
| 9:2-07a | 3-09 | AAAAbAbbb | 7436529 | 0 | 9 | 1 | 9-2-05 |
| 9:2-07b | 5n11 | AABaBcBcdCd | 2916 | 2 | 9 | 1 | 9-2-05+- |
| 9:2-08a | 3-09 | AAAAbbAbb | 7527429 | 0 | 2 | 4 | 9-2-14+- |
| 9:2-08b | 4n10 | AAAAbACBBC | 237267 | 1 | 2 | 4 | 9-2-14 |

Table 2. Cont'd.

| | | | | | | | |
|---|---|---|---|---|---|---|---|
| 9:2-09a | 3-09 | AAAbAAAbb | 7527429 | 0 | 2 | 2 | 9-2-23 |
| 9:2-09b | 5n13 | AABaBCbbdCBdC | 2736 | 3 | 7 | 2 | 9-2-23+- |
| 9:2-10a | 3-09 | AAAbAAbAb | 6709239 | 0 | 1 | 2 | 9-2-35 |
| 9:2-10b | 4n10 | AAbAbAbcBc | 131778 | 1 | 8 | 2 | 9-2-35+- |
| 9:2-11a | 3-09 | AAAbAAbbb | 7518429 | 0 | 1 | 2 | 9-2-21+- |
| 9:2-11b | 4n10 | AABaBcBBBc | 139068 | 1 | 8 | 2 | 9-2-21 |
| 9:2-12a | 3-09 | AAAbAbAbb | 6691239 | 0 | 8 | 1 | 9-2-34 |
| 9:2-12b | 4n10 | AABacBBcBc | 133578 | 1 | 1 | 1 | 9-2-34+- |
| 9:2-13 | 3-09 | AAAbAbbAb | 6700239 | 0 | 9 | 2 | 9-2-37 |
| 9:2-14a | 3-09 | AAAbbAAbb | 7600329 | 0 | 2 | 4 | 9-2-29+- |
| 9:2-14b | 4n10 | AAAbAACBBC | 229977 | 1 | 2 | 4 | 9-2-29 |
| 9:2-15a | 3-09 | AAbAAbAAb | 6800139 | 0 | 3 | 4 | 9-2-40+- |
| 9:2-15b | 5n13 | AABaCbACBBDcD | 5373 | 3 | 3 | 4 | 9-2-40 |
| 9:2-16a | 3-09 | AAbAAbAbb | 6782139 | 0 | 1 | 2 | 9-2-39 |
| 9:2-16b | 4n10 | AABBcBaBcB | 124488 | 1 | 8 | 2 | 9-2-39+- |
| 9:2-17a | 3-09 | AAbAbAAbb | 6773139 | 0 | 9 | 2 | 9-2-41 |
| 9:2-17b | 5n11 | AbACBBdcBcd | 37179 | 1 | 9 | 2 | 9-2-41+- |
| 9:2-18a | 3-09 | AAbAbAbAb | 5954949 | 0 | 8 | 2 | 9-2-42 |
| 9:2-18b | 4n10 | AABcBaBcBc | 118998 | 1 | 1 | 2 | 9-2-42+- |
| 9:2-19a | 5-09 | AAbACbdCd | 43569 | 0 | 8 | 2 | 9-2-11+- |
| 9:2-19b | 4n10 | AAABaBBcBc | 131778 | 1 | 8 | 2 | 9-2-11 |
| 9:2-20a | 5-09 | AbAbCbdCd | 45369 | 0 | 1 | 1 | 9-2-12+- |
| 9:2-20b | 4n10 | AAABaBcBcc | 133578 | 1 | 1 | 1 | 9-2-12 |
| 9:2-21 | 5-09 | AbACbbdCd | 44469 | 0 | 9 | 2 | 9-2-25 |
| 9:2-22a | 4n10 | AAAAABacBc | 146358 | 1 | 8 | 3 | 9-2-03 |
| 9:2-22b | 5n11 | AAABacBcdCd | 3006 | 2 | 1 | 3 | 9-2-03+- |
| 9:2-23 | 4n10 | AAAAbAbcBc | 154548 | 1 | 9 | 2 | 9-2-04 |
| 9:2-24a | 4n10 | AAAABaBcBc | 130878 | 1 | 7 | 2 | 9-2-07 |
| 9:2-24b | 4n10 | AAbAbbbcBc | 141768 | 1 | 2 | 2 | 9-2-07+- |
| 9:2-25 | 4n10 | AAAABacBBc | 139968 | 1 | 9 | 3 | 9-2-27 |
| 9:2-26a | 4n10 | AAAABacBcc | 149058 | 1 | 2 | 2 | 9-2-06+- |
| 9:2-26b | 4n10 | AABaBBBcBc | 138168 | 1 | 7 | 2 | 9-2-06 |
| 9:2-27 | 4n10 | AAAbAbbcBc | 139968 | 1 | 9 | 2 | 9-2-09 |
| 9:2-28 | 4n10 | AAABaBcBBc | 132678 | 1 | 9 | 2 | 9-2-18 |
| 9:2-29a | 4n10 | AAAbAbcBcc | 148158 | 1 | 1 | 2 | 9-2-08 |
| 9:2-29b | 5n11 | AABaBBCbdCd | 2826 | 2 | 8 | 2 | 9-2-08+- |
| 9:2-30a | 4n10 | AAAbACBBBC | 230877 | 1 | 3 | 4 | 9-2-22+- |
| 9:2-30b | 4n10 | AABaBcBBcc | 142668 | 1 | 3 | 4 | 9-2-22 |
| 9:2-31a | 4n10 | AAABacBBcc | 141768 | 1 | 2 | 3 | 9-2-28 |
| 9:2-31b | 5n11 | AABaBAcBDcD | 42669 | 1 | 7 | 3 | 9-2-28+- |
| 9:2-32 | 4n10 | AAABacBcBc | 125388 | 1 | 9 | 2 | 9-2-33 |
| 9:2-33a | 4n10 | AAABaccBcc | 149058 | 1 | 2 | 3 | 9-2-17+- |
| 9:2-33b | 5n11 | AAABacBDCCD | 3825 | 2 | 2 | 3 | 9-2-17 |
| 9:2-34a | 4n10 | AAABBaBcBc | 123588 | 1 | 7 | 3 | 9-2-26 |
| 9:2-34b | 5n11 | ABaBAcBcDcD | 46269 | 1 | 2 | 3 | 9-2-26+- |
| 9:2-35a | 4n10 | AAABcBaBcB | 116298 | 1 | 7 | 2 | 9-2-38+- |
| 9:2-35b | 4n10 | AAbCbAbccb | 127188 | 1 | 2 | 2 | 9-2-38 |
| 9:2-36a | 4n10 | AAbAAbbcBc | 141768 | 1 | 2 | 3 | 9-2-16+- |
| 9:2-36b | 5n11 | AABaBcBDCCD | 4554 | 2 | 2 | 3 | 9-2-16 |
| 9:2-37a | 4n10 | AAbAACBBBC | 222687 | 1 | 2 | 3 | 9-2-30+- |
| 9:2-37b | 5n11 | AAbACBBdccd | 35379 | 1 | 7 | 3 | 9-2-30 |
| 9:2-38 | 4n10 | AABaBBcBBc | 139968 | 1 | 9 | 3 | 9-2-15 |
| 9:2-39 | 4n10 | AAbAbcBAcb | 132678 | 1 | 9 | 3 | 9-2-36 |

Table 2. Cont'd.

| | | | | | | | |
|---|---|---|---|---|---|---|---|
| 9:2-40 | 4n10 | AABaBcBcBc | 125388 | 1 | 9 | 2 | 9-2-32 |
| 9:2-41a | 4n10 | AAbACBBBCC | 223587 | 1 | 3 | 3 | 9-2-24+- |
| 9:2-41b | 5n11 | AbaCbaCbdCd | 47169 | 1 | 3 | 3 | 9-2-24 |
| 9:2-42 | 5n11 | AAABaBCbdCd | 2187 | 2 | 9 | 2 | 9-2-10 |
| 9:2n43a | 3n09 | AAAAAbaab | 90009 | 1 | 2 | 4 | 9-2-43 |
| 9:2n43b | 3n09 | AAAAABAAB | 9000009 | 0 | 2 | 4 | 9-2-43+- |
| 9:2n44 | 3n09 | AAAAABaaB | 73629 | 1 | 9 | 2 | 9-2-44 |
| 9:2n45a | 3n09 | AAAAbaaab | 81819 | 1 | 1 | 2 | 9-2-50 |
| 9:2n45b | 4n10 | AAABaaBCbC | 639 | 2 | 8 | 2 | 9-2-50+- |
| 9:2n46a | 3n09 | AAAABAAAB | 9009009 | 0 | 3 | 4 | 9-2-49 |
| 9:2n46b | 4n10 | AAAbaabcBc | 999 | 2 | 3 | 4 | 9-2-49+- |
| 9:2n47a | 3n09 | AAAABaaaB | 83619 | 1 | 3 | 4 | 9-2-51 |
| 9:2n47b | 4n10 | AAABAABCbC | 171828 | 1 | 3 | 4 | 9-2-51+- |
| 9:2n48a | 3n09 | AAAABaaBB | 65439 | 1 | 8 | 2 | 9-2-52 |
| 9:2n48b | 4n10 | AABaaBAcBc | 2277 | 2 | 1 | 2 | 9-2-52+- |
| 9:2n49 | 3n09 | AAAABaBaB | 66339 | 1 | 9 | 2 | 9-2-55 |
| 9:2n50a | 3n09 | AAAbaaabb | 74529 | 1 | 1 | 1 | 9-2-54+- |
| 9:2n50b | 4n10 | AAbAbcBaBc | 1368 | 2 | 8 | 1 | 9-2-54 |
| 9:2n51a | 3n09 | AAABAAABB | 9018009 | 0 | 4 | 4 | 9-2-53 |
| 9:2n51b | 4n10 | ABaBACBaBC | 1089 | 2 | 4 | 4 | 9-2-53+- |
| 9:2n52a | 3n09 | AAABaaBaB | 75429 | 1 | 2 | 2 | 9-2-57+- |
| 9:2n52b | 4n10 | AAAbaaCBBC | 64539 | 1 | 7 | 2 | 9-2-57 |
| 9:2n53a | 3n09 | AAABaBaBB | 57249 | 1 | 7 | 2 | 9-2-58+- |
| 9:2n53b | 4n10 | AAABaacbbc | 68139 | 1 | 2 | 2 | 9-2-58 |
| 9:2n54a | 3n09 | AAABBAABB | 8927109 | 0 | 2 | 4 | 9-2-59+- |
| 9:2n54b | 4n10 | AAABAAcbbc | 97299 | 1 | 2 | 4 | 9-2-59 |
| 9:2n55 | 3n09 | AABaBAABB | 99099 | 1 | 4 | 4 | 9-2-61 |
| 9:2n56 | 4n10 | AAABaBCbbC | 729 | 2 | 9 | 1 | 9-2-47 |
| 9:2n57a | 4n10 | AAABaBCBBC | 163638 | 1 | 2 | 3 | 9-2-48+- |
| 9:2n57b | 4n10 | AAABaBcbbc | 2367 | 2 | 2 | 3 | 9-2-48 |
| 9:2n58 | 4n10 | AAAbACbaCB | 66339 | 1 | 9 | 2 | 9-2-56 |
| 9:2n59a | 4n10 | AAAbaCbaCb | 1638 | 2 | 2 | 3 | 9-2-45+- |
| 9:2n59b | 4n10 | AAABACbACB | 170928 | 1 | 2 | 3 | 9-2-45 |
| 9:2n60 | 4n10 | AAbaabbcBc | 1458 | 2 | 9 | 2 | 9-2-46 |
| 9:2n61a | 4n10 | AABAAcbbbc | 90009 | 1 | 2 | 3 | 9-2-60+- |
| 9:2n61b | 4n10 | AABaBACBBC | 97299 | 1 | 2 | 3 | 9-2-60 |
| 9:3-01a | 4-09 | AAAbACbbC | 588951 | 0 | 8 | 4 | 9-3-06 |
| 9:3-01b | 4n11 | AAABBaBBCbC | 30294 | 2 | 5 | 4 | 9-3-06++- |
| 9:3-01c | 5n12 | AABaBBCbdccd | 10692 | 2 | 6 | 4 | 9-3-06+-+ |
| 9:3-02a | 4-09 | AAbAAbCbC | 605151 | 0 | 1 | 3 | 9-3-02+-- |
| 9:3-02b | 4-09 | AAbAbCbbC | 605151 | 0 | 1 | 3 | 9-3-02++- |
| 9:3-02c | 5n10 | AbAbCbdccd | 18873 | 1 | 8 | 3 | 9-3-02 |
| 9:3-03a | 4-09 | AAbAACbbC | 597051 | 0 | 9 | 3 | 9-3-08++- |
| 9:3-03b | 5n10 | AbACbbdccd | 19683 | 1 | 9 | 3 | 9-3-08 |
| 9:3-04a | 4-09 | AAbAACbCC | 588951 | 0 | 8 | 3 | 9-3-01+-+ |
| 9:3-04b | 4-09 | AbAbbCbbC | 588951 | 0 | 8 | 3 | 9-3-01 |
| 9:3-04c | 5n10 | AAbACbdccd | 20493 | 1 | 1 | 3 | 9-3-01++- |
| 9:3-05a | 4-09 | AAbACbACb | 523341 | 0 | 8 | 4 | 9-3-11 |
| 9:3-05b | 4-09 | AAbCbAbCb | 621351 | 0 | 3 | 4 | 9-3-11+-- |
| 9:3-05c | 4n11 | AAABcBaBCCB | 36855 | 2 | 5 | 4 | 9-3-11++- |
| 9:3-06a | 4-09 | AAbACbCbC | 531441 | 0 | 9 | 2 | 9-3-10++- |
| 9:3-06b | 4-09 | AbAbAbCbC | 531441 | 0 | 9 | 2 | 9-3-10 |
| 9:3-07a | 4-09 | AAbbACbbC | 605151 | 0 | 1 | 4 | 9-3-04++- |

Table 2. Cont'd.

| | | | | | | | |
|---|---|---|---|---|---|---|---|
| 9:3-07b | 4n11 | AAABaaCbbCb | 605151 | 1 | 1 | 4 | 9-3-04+-+ |
| 9:3-07c | 5n12 | AABaBAcBDCCD | 28674 | 2 | 3 | 4 | 9-3-04 |
| 9:3-08a | 4-09 | AbAbACbbC | 531441 | 0 | 9 | 2 | 9-3-12+-+ |
| 9:3-08b | 4n11 | AAAbbacBcBc | 531441 | 1 | 9 | 2 | 9-3-12 |
| 9:3-09a | 4-09 | AbACbbbbC | 572751 | 0 | 6 | 4 | 9-3-05+-- |
| 9:3-09b | 3n10 | AAAAABBaBB | 1562652 | 1 | 5 | 4 | 9-3-05 |
| 9:3-09c | 5n10 | AABacBDccD | 12312 | 1 | 8 | 4 | 9-3-05++- |
| 9:3-10 | 4-09 | AbbAbCbbC | 597051 | 0 | 9 | 3 | 9-3-09 |
| 9:3-11a | 4-09 | AbbACbbbC | 588951 | 0 | 8 | 4 | 9-3-03+-+ |
| 9:3-11b | 4n11 | AAAAbbacBBc | 588951 | 1 | 8 | 4 | 9-3-03 |
| 9:3-11c | 4n11 | AAAABBaCBBC | 1480842 | 1 | 3 | 4 | 9-3-03+-- |
| 9:3-12a | 5n10 | AbAbcBDccD | 12312 | 1 | 8 | 3 | 9-3-07 |
| 9:3-12b | 4n11 | AAABaBCCbCC | 22113 | 2 | 3 | 3 | 9-3-07++- |
| 9:3n13a | 4n09 | AAbaabCbC | 5751 | 1 | 8 | 3 | 9-3-14 |
| 9:3n13b | 4n09 | AABAAcBcc | 736371 | 0 | 1 | 3 | 9-3-14++- |
| 9:3n13c | 4n09 | AAbAbCBBC | 5751 | 1 | 8 | 3 | 9-3-14+-- |
| 9:3n14a | 4n09 | AABAABcBc | 720171 | 0 | 8 | 3 | 9-3-13+-- |
| 9:3n14b | 4n09 | AAbaaCbCC | 7371 | 1 | 1 | 3 | 9-3-13++- |
| 9:3n14c | 4n09 | AAbacBacb | 7371 | 1 | 1 | 3 | 9-3-13 |
| 9:3n15 | 4n09 | AAbaaCbbC | 6561 | 1 | 9 | 2 | 9-3-18 |
| 9:3n16a | 4n09 | AABAAcBBc | 728271 | 0 | 9 | 3 | 9-3-19++- |
| 9:3n16b | 4n09 | ABaBACbbC | 6561 | 1 | 9 | 3 | 9-3-19 |
| 9:3n17a | 4n09 | AAbbACBBC | 13932 | 1 | 1 | 4 | 9-3-16 |
| 9:3n17b | 4n09 | AAbbAcbbc | 752571 | 0 | 3 | 4 | 9-3-16+-+ |
| 9:3n17c | 4n09 | AbbACBBBC | 13932 | 1 | 1 | 4 | 9-3-16+-- |
| 9:3n18a | 4n09 | AAbbAcBBc | 81 | 2 | 1 | 4 | 9-3-17+-+ |
| 9:3n18b | 4n09 | AABBACBBC | 883791 | 0 | 3 | 4 | 9-3-17++- |
| 9:3n18c | 4n11 | AAAAbbaCbbC | 81 | 3 | 1 | 4 | 9-3-17 |
| 9:3n19a | 4n09 | AABBAcBBc | 736371 | 0 | 1 | 4 | 9-3-15+-+ |
| 9:3n19b | 4n09 | AbbAcbbbc | 736371 | 0 | 1 | 4 | 9-3-15++- |
| 9:3n19c | 4n09 | ABBACbbbC | 15552 | 1 | 3 | 4 | 9-3-15 |
| 9:3n20a | 4n09 | AbAbACBBC | 13122 | 1 | 9 | 3 | 9-3-20 |
| 9:3n20b | 4n09 | AbAbAcbbc | 662661 | 0 | 9 | 3 | 9-3-20++- |
| 9:3n21a | 4n09 | AbAbAcBBc | 0 | 2 | 9 | 1 | 9-3-21+-+ |
| 9:3n21b | 4n11 | AAAbbaCbCbC | 0 | 3 | 9 | 1 | 9-3-21 |
| 9:4-01a | 4n10 | AAABBacBBc | 132678 | 1 | 2 | 4 | 9-4-01++-- |
| 9:4-01b | 5n11 | ABaBAcBDccD | 44469 | 1 | 7 | 4 | 9-4-01 |

Table 3. Minimum Braids for "Trivial" Links

| Tag # | S-Cr | Minimum Braid | AP(10) | z | d | u | Description |
|---|---|---|---|---|---|---|---|
| 4:3.01 | 3-04 | AAbb | 81 | 0 | 1 | 2 | L2+1 |
| 5:2.01 | 3-05 | AAAbb | 819 | 0 | 1 | 2 | K3+1 |
| 6:1c01a | 3-06 | AAAbbb | 8281 | 0 | 1 | 2 | Square |
| 6:1c01b | 3-06 | AAABBB | 8281 | 0 | 1 | 2 | Granny |
| 6:2.01 | 4-06 | AAbCbC | 639 | 0 | 8 | 2 | 1+K4 |
| 6:3.01a | 3-06 | AAAAbb | 8181 | 0 | 2 | 3 | L4a+1 |
| 6:3.01b | 4n07 | AABaBCC | 162 | 1 | 2 | 3 | L4b+1 |
| 6:4.01 | 4-06 | AAbbCC | 729 | 0 | 1 | 3 | L2+1+1 |
| 6:4.02 | 4-06 | AAbCCb | 729 | 0 | 1 | 3 | 1+L2+1 |
| 7:1c01 | 4-07 | AAAbCbC | 6461 | 0 | 8 | 2 | K3#K4 |
| 7:2c01 | 3-07 | AAAAbbb | 82719 | 0 | 2 | 3 | L4a#K3 |
| 7:2c02 | 4n08 | AAAbbcBc | 1638 | 1 | 2 | 3 | K3#L4b |
| 7:2.01 | 3-07 | AAAAAbb | 81819 | 0 | 1 | 3 | K5.1+1 |
| 7:2.02 | 4n08 | AAABaBCC | 1548 | 1 | 1 | 2 | K5.2+1 |
| 7:3.01 | 4-07 | AAAbbCC | 7371 | 0 | 1 | 3 | K3+1+1 |
| 7:3.02 | 4-07 | AAAbCCb | 7371 | 0 | 1 | 3 | 1+K3+1 |
| 7:3.03 | 4-07 | AAbAbCC | 6561 | 0 | 9 | 2 | L5+1 |
| 8:1c01 | 3-08 | AAAAAbbb | 827281 | 0 | 1 | 3 | K5.1#K3 |
| 8:1c02 | 5-08 | AbAbCdCd | 5041 | 0 | 1 | 2 | K4#K4 |
| 8:1c03 | 4n09 | AAABaBCCC | 15652 | 1 | 1 | 2 | K5.2#K3 |
| 8:2c01 | 4-08 | AAAAbCbC | 64539 | 0 | 7 | 3 | L4a#K4 |
| 8:2c02 | 4-08 | AAAbbCbC | 66339 | 0 | 9 | 2 | K3#K5.1 |
| 8:2c03 | 5n09 | AABaBCdCd | 1278 | 1 | 7 | 3 | L4b#K4 |
| 8:2.01 | 4-08 | AAAbAbCC | 65439 | 0 | 8 | 2 | K6.1+1 |
| 8:2.02 | 4-08 | AAAbbbCC | 74529 | 0 | 1 | 3 | K3#K3+1 |
| 8:2.03 | 4-08 | AAbAbbCC | 67239 | 0 | 1 | 2 | K6.2+1 |
| 8:2.04 | 5n09 | AABacBcdd | 1368 | 1 | 8 | 2 | K6.3+1 |
| 8:2_01 | 4-08 | AAAbbCCC | 74529 | 0 | 1 | 3 | K3+K3 |
| 8:3c01 | 3-08 | AAAAbbbb | 826281 | 0 | 4 | 4 | L4a#L4a |
| 8:3c02 | 4n09 | AAAAbbcBc | 16362 | 1 | 4 | 4 | L4a#L4b |
| 8:3c03 | 5n10 | AABaBCCDcD | 324 | 2 | 4 | 4 | L4b#L4b |
| 8:3.01a | 3-08 | AAAAAAbb | 818181 | 0 | 3 | 4 | L6.1a+1 |
| 8:3.01b | 5n10 | AABaBCbCdd | 243 | 2 | 3 | 4 | L6.1b+1 |
| 8:3.02 | 5-08 | AAbbCdCd | 5751 | 0 | 8 | 3 | 1+1+K4 |
| 8:3.03a | 5-08 | AAbCbdCd | 4941 | 0 | 7 | 3 | 1+L6.2a |
| 8:3.03b | 4n09 | AAABaBBCC | 14742 | 1 | 2 | 3 | L6.2b+1 |
| 8:3.04 | 4n09 | AAAABaBCC | 15552 | 1 | 3 | 4 | L6.3+1 |
| 8:4.01a | 4-08 | AAAAbbCC | 73629 | 0 | 2 | 4 | L4a+1+1 |
| 8:4.01b | 5n09 | AABaBCCdd | 1458 | 1 | 2 | 4 | L4b+1+1 |
| 8:4.02a | 4-08 | AAAAbCCb | 73629 | 0 | 2 | 4 | 1+L4a+1 |
| 8:4.02b | 5n09 | AABaBCddC | 1458 | 1 | 2 | 4 | 1+L4b+1 |
| 8:4.03a | 4-08 | AAbAAbCC | 66339 | 0 | 1 | 4 | 6:3-01a+1 |
| 8:4.03b | 5n09 | AAbCbdccd | 2187 | 1 | 3 | 4 | 1+6:3-01b |
| 8:4.04a | 4-08 | AAbbbbCC | 73629 | 0 | 2 | 4 | 1,L4a+1 |
| 8:4.04b | 5n09 | AAbbcBddc | 1458 | 1 | 2 | 4 | 1,L4b+1 |
| 8:4.05 | 4-08 | AAbCbCbC | 59049 | 0 | 9 | 3 | 1+6:3-02 |
| 8:4:01a | 4n08 | AAbaabCC | 729 | 1 | 1 | 4 | 6:3n03a+1 |
| 8:4:01b | 4n08 | AABAABCC | 80919 | 0 | 3 | 4 | 6:3n03b+1 |
| 8:5.01 | 5-08 | AAbbCCdd | 6561 | 0 | 1 | 4 | L2+1+1+1 |
| 8:5.02 | 5-08 | AAbbCddC | 6561 | 0 | 1 | 4 | 1+1+L2+1 |
| 8:5.03 | 5-08 | AAbCddCb | 6561 | 0 | 1 | 4 | 1+(L2+1)+1 |

Table 3. Cont'd.

| | | | | | | | |
|---|---|---|---|---|---|---|---|
| 9:1c01 | 4-09 | AAAAAbCbC | 645461 | 0 | 8 | 3 | K5.1#K4 |
| 9:1c02 | 4-09 | AAAbAbCCC | 661661 | 0 | 8 | 2 | K6.1#K3 |
| 9:1c03 | 4-09 | AAAbbbCCC | 753571 | 0 | 1 | 3 | K3#K3#K3 |
| 9:1c04 | 4-09 | AAAbbCbCC | 679861 | 0 | 1 | 2 | K3#K6.2 |
| 9:1c05 | 5n10 | AAABaBCdCd | 12212 | 1 | 8 | 2 | K5.2#K4 |
| 9:1c06 | 5n10 | AAAbbcBDcD | 13832 | 1 | 8 | 2 | K3#K6.3 |
| 10:1c01 | 3-10 | AAAAAAAbbb | 82727281 | 0 | 1 | 4 | K7.1#K3 |
| 10:1c02 | 3-10 | AAAAAbbbbb | 82646281 | 0 | 1 | 4 | K5.1#K5.1 |
| 10:1c03 | 5-10 | AAAbAbCdCd | 516241 | 0 | 1 | 2 | K6.1#K4 |
| 10:1c04 | 5-10 | AAAbbbCdCd | 587951 | 0 | 8 | 3 | K3#K3#K4 |
| 10:1c05 | 5-10 | AAAbbCbdCd | 514241 | 0 | 8 | 2 | K3#K7.2 |
| 10:1c06 | 5-10 | AAAbCbCdCd | 532441 | 0 | 1 | 2 | K3#K7.3 |
| 10:1c07 | 5-10 | AAbAbbCdCd | 530441 | 0 | 8 | 2 | K6.2#K4 |
| 10:1c08 | 4n11 | AAAAABaBCCC | 1571752 | 1 | 1 | 3 | K7.4#K3 |
| 10:1c09 | 4n11 | AAAAAbbbcBc | 1563652 | 1 | 1 | 3 | K5.1#K5.2 |
| 10:1c10 | 4n11 | AAAABaBBCCC | 1498042 | 1 | 1 | 3 | K7.5#K3 |
| 10:1c11 | 6n11 | AABacBcdEdE | 10792 | 1 | 1 | 2 | K6.3#K4 |
| 10:1c12 | 5n12 | AAABaBCbCddd | 23023 | 2 | 1 | 2 | K7.6#K3 |
| 10:1c13 | 5n12 | AAABaBCCCDcD | 29584 | 2 | 1 | 2 | K5.2#K5.2 |
| 10:1c14 | 5n12 | AAAbbcBccdCd | 30394 | 2 | 1 | 3 | K3#K7.7 |

Table 5 summarizes the differences between the minimum braids and the Chalcraft braids. In most cases, the Chalcraft braids have more crossings and in 15 cases more strands than the minimum braids.

Table 5. Differences for Chalcraft's braids

| | | Differences | | | | | | | | | | | |
|---|---|---|---|---|---|---|---|---|---|---|---|---|---|
| | | Strands | | Crossings | | | | | | | | | |
| Strands | Knots | 0 | 1 | 0 | 1 | 2 | 3 | 4 | 5 | 6 | 7 | 8 | 9 |
| 2 | 4 | 4 | - | 4 | - | - | - | - | - | - | - | - | - |
| 3 | 58 | 55 | 3 | 55 | - | - | 3 | - | - | - | - | - | - |
| 4 | 116 | 105 | 11 | 19 | - | 86 | 3 | - | 4 | - | 3 | - | 1 |
| 5 | 66 | 65 | 1 | 3 | - | 11 | - | 45 | - | 6 | - | - | 1 |
| 6 | 5 | 5 | - | - | - | - | - | - | - | - | - | 5 | - |
| Totals | 249 | 234 | 15 | 81 | 0 | 97 | 6 | 45 | 4 | 6 | 3 | 5 | 2 |

These two alternative braid sets are useful when testing ways to reduce nonminimum braids to their minimum by using Reidemeister and Markov moves. A forthcoming paper will describe how these moves can reduce the Jones braids to minimum braids. [11]

## 4. Graph Trees

The first application of minimum braids is a mapping of a set of links to a set of trees in graph theory so as to assign a unique minimum braid to each tree. These braids are then used to propose a conjecture about the number of trees with alternating minimum braids and the number of rational links in knot theory. The Alexander polynomials of these braids have some interesting properties concerning the weighted sums of their coefficients. In the section on periodic tables, these weighted sums are used to classify families of links.

The links that are mapped to trees consist of simple unknots or loops that are connected by two alternating crossings and are acyclic. If the number of components is equal to n, then the number of crossings is equal to 2*(n-1). The mapping projects each loop to a vertex in the tree and each *pair* of link crossings to a *single* edge that connects the corresponding pair of vertices. With this mapping, the number of vertices in the graph is equal to one plus the number of edges. This is one of the definitions of a tree [12]. The usual mapping between links and planar graphs makes the number of link crossings equal to the number of edges in the graph.

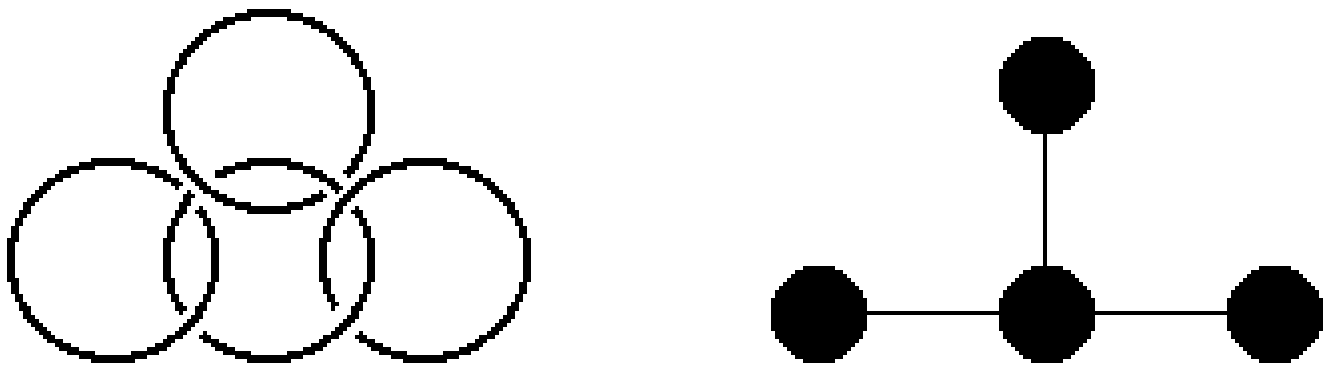

Figure 1. Mapping Links to Trees

In Figure 1, the number of components in the link and the number of vertices in the tree is equal to four. The number of link crossings is six and the number of edges in the tree is three. In Table 3, this link has the Tag # of 6:4.02 and the minimum braid AAbCCb.

For all trees through six vertices, the minimum braid is always alternating. The first nonalternating minimum braid is for one of the trees with seven vertices. Figure 2 plots the 11 trees with seven vertices.

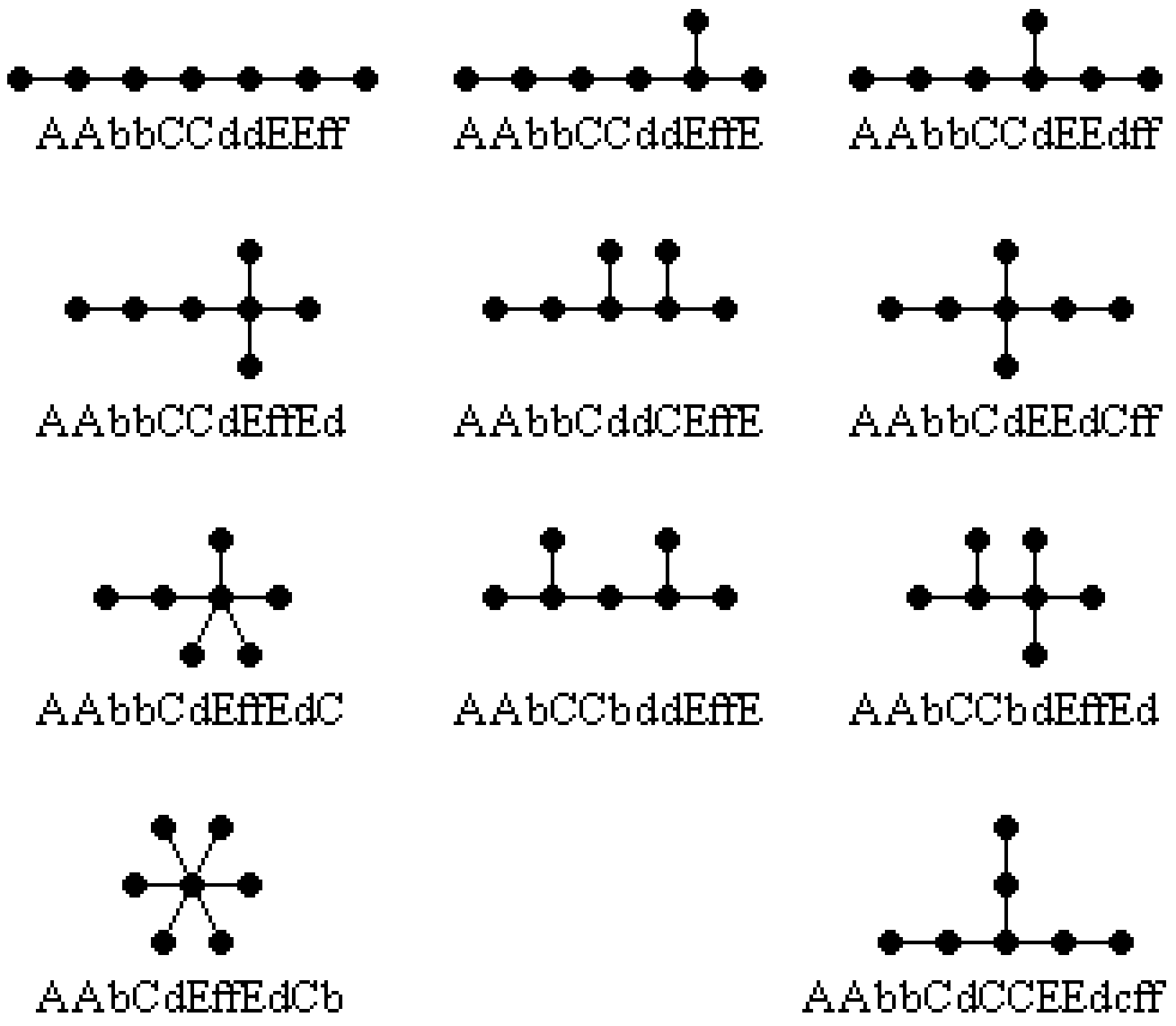

Figure 2. Minimum Braids of Trees with Seven Vertices

The trees are ordered according to their minimum braids, which are listed beneath each tree. The first ten trees have alternating minimum braids with 12 braid crossings. The last tree has a 14 crossing minimum braid that is nonalternating. This tree is the first tree with a vertex that has more than two branches of lengths greater than one. An examination of trees through ten crossing found that this type of tree always has a nonalternating minimum braid.

For trees with 8-10 vertices, the number of nonalternating minimum braids are 3, 11, and 34, respectively. Starting with a tree with one vertex, the number of trees with alternating minimum braids is 1, 1, 1, 2, 3, 6, 10, 20, 36, and 72. As has been noted by Jablan, this is the integer sequence for the number of rational knots and links. [13]

**Conjecture 1.** The number of graph trees of n vertices with alternating minimum braids is equal to the number of rational knots and links with n crossings.

When counting the number of rational knots and links, the orientations are not differentiated. A future paper will present a simple one-to-one mapping of trees with alternating minimum braids and rational knots and links.

Another interesting property of minimum braids of trees concerns the coefficients of the Alexander polynomials. Since the Alexander polynomial (and all of the other existing polynomials) can not distinguish between these "trivial" links that are mapped to trees, all of the trees have the same coefficients for their polynomials. The Alexander polynomial coefficients are equal to the Mth row of a negative Pascal triangle, and therefore have some well known properties. [14] If the weighted sums of the coefficients (c(i)) of the Alexander polynomial is defined as

$W[M,j] = \sum i^j c(i)$, $i = 1..M$, then
$W[M,j] = 0$, $0 <= j < (M-1)$, and
$W[M,j] = s(M,1)\ S(j+1,M)$, $j >= (M-1)$,

where s(n,m) represents Stirling numbers of the first kind (Riordan notation) and S(n,m) represents Stirling numbers of the second kind. In the section on periodic tables, weighted sums of Alexander polynomials are used to classify families of links.

## 5. Amphicheiral Links

A second application of minimum braids is the study of symmetries of knots and orientated links. In general, braids can have symmetry because they can always be read from top to bottom and then rotated 180 degrees and again read from top to bottom. Sometimes symmetry can be found in a minimum braid, while at other times, the minimum braid has to be transformed using Reidemeister-Markov moves to get it into a nonminimum but symmetric form.

The simplest symmetry to find this way is amphicheirality without regards to invertibility. Links with this symmetry are equal to their mirror images and are unchanged when all of the crossings are reversed or changed to the opposite type of crossings. If a braid that has been rotated is equal to the reverse of the original braid, then the corresponding knot or orientated link is amphicheiral. This is a sufficient condition, but not a necessary one. Such a braid is said to be a reverse rotated palindrome (RRP).

There are 20 amphicheiral knots with ten or less crossings and four amphicheiral links with eight or less crossings. The minimum braids of 13 of the knots and two of the braids are RRPs and the amphicheirality can be read directly from the braids. These knots are 4:1-01, 6:1-02, 8:1-05, 8:1-07, 8:1-08, 8:1-09, 10:1-009, 10:1-017, 10:1-020, 10:1-023, 10:1-028, 10:1-031, and 10:1-033. The links are 6:3-02 and 8:3-05a. The remaining braids have to be transformed to get them into palindromic form. Table 6 lists these knots and links, their minimum braids, and one of their reverse rotated palindromes.

Table 6. Reverse Rotated Palindromes

| Tag # | Minimum Braid | Palindrome |
|---|---|---|
| 8:1-18 | AABacBcdCd | AABacBDcdd |
| 8:3-04a | AAbAbbAb | AbbAbAAb |
| 8:3-05c | AbACBBdccd | AbAccBBdCd |
| 10:1-021 | AAbAbAbbAb | AbAAbAbbAb |
| 10:1-022 | AAbAbbAAbb | AAbbAbAAbb |
| 10:1-093 | AAABacBcdCdd | AAABacBDcddd |
| 10:1-103 | AABaBcBccdCd | AABaBcBcDcdd |
| 10:1-109 | AAbACBBdcccd | AAbAccBBdCdd |
| 10:1-115 | AbACBBdcBccd | AbAccBcBBdCd |

Of the nine minimum braids listed in Table 6, the only difficult one to transform into a RRP is 10:1-109. This transformation requires Markov and Reidemeister moves

In order for a RRP to exist there must be an even number of braid crossings in the minimum braid. There can be an odd or an even number of crossings in the knot or link. In 1998, Hoste, Thistlethwaite, and Weeks [15] published the first amphicheiral knot with an odd number of crossings, and thereby disproved a hundred year old conjecture by Tait. [16] Using the symmetry conditions of RRP, it was easy to find a five strand, 16 crossing braid for this 15 crossing knot. The RRP is $ABaBC^nBAdcb^ncDcd$, where n is equal to two. By increasing the value for n, it is possible to generate a infinite family of odd-crossing knots that are amphicheiral.

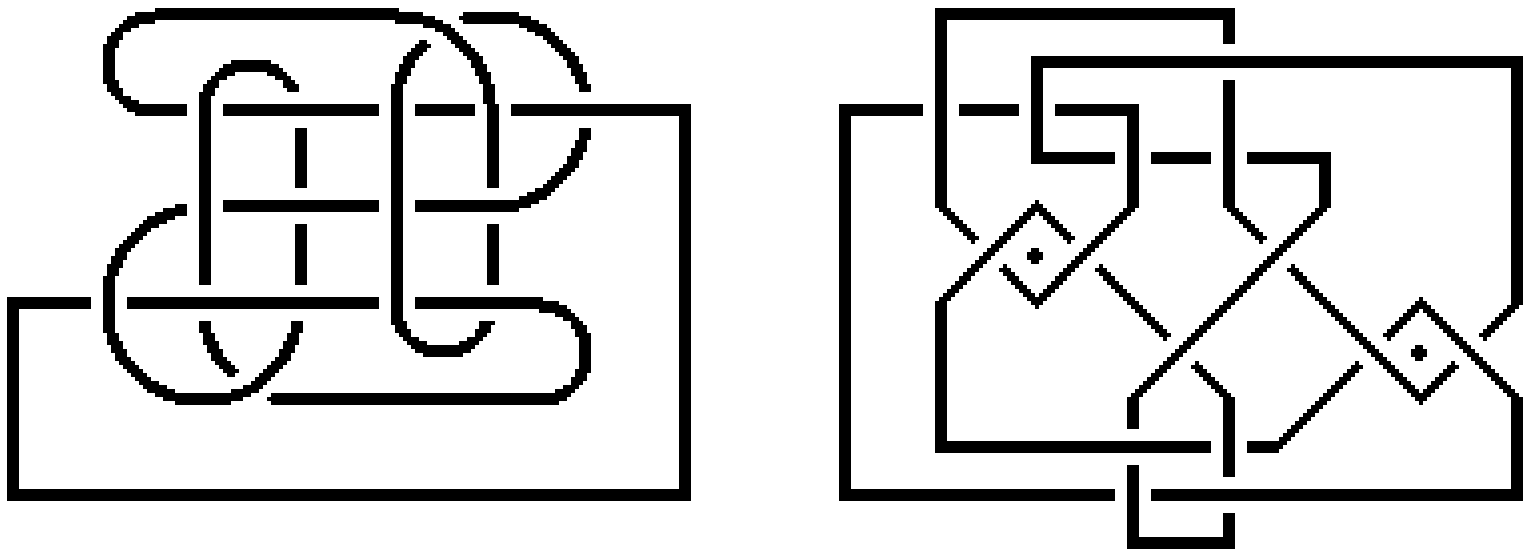

Figure 3. The French Horn Knot

Figure 3 plots a 16 crossing projection of this knot that shows the symmetry. On the right side of this figure is a 15 crossing projection where dots have been places to show where additional crossings are generated when the value of n is increased. Due to the symmetry of the 16 crossing projection, this knot can be called the French horn knot.

Another form of symmetry of braid is a rotated palindrome that is not reversed. The braid reads exactly the same before and after the rotation. This type of braid can be associated with other types of symmetry and is a subject for future studies. An example of a rotated palindrome is ABBAcBaDcbAdCbDCCD. This is a 16 crossing knot and the first knot to have $D_9$ symmetry [15].

## 6. Unknotting Numbers

A third application of minimum braids involves different approaches to the calculation of unknotting numbers. An unknotting number is the minimum number of crossings that must be reversed or switched in order to reduce the original knot to the unknot. There are two different ways to specify unknotting problems. The classical approach allows Reidemeister-Markov moves after each reversal of a crossing. A special type of unknotting number allows these moves only after *all* reversals have been made. The most famous example of this approach is the 514 pretzel knot which has a classic unknotting number of two. Bleiler and Nalsanishi proved independently in 1983 that the special unknotting number for any ten crossing projection of this knot is three.

Using minimum braids it is easy to calculate now many *braid* crossings must be changed before the braid becomes the unknot. This invariant is the minimum braid unknotting number. It can be determined by systematically reversing every possible combination of braid crossings. It is not known if this is a new invariant or just a simple way to calculate the classical unknotting number.

In Table 1 the 514 pretzel knot is listed as knot 10:1-035 and has the minimum braid AAAAAbAbcBc. Reverse two of the braid crossings to get an 11 crossing projection (AAAAAbabcbc) of the unknot. Figure 4 plots the 10 and an 11 crossing projection of this knot. The two crossings that can be reversed in the 11 crossing projection are circled.

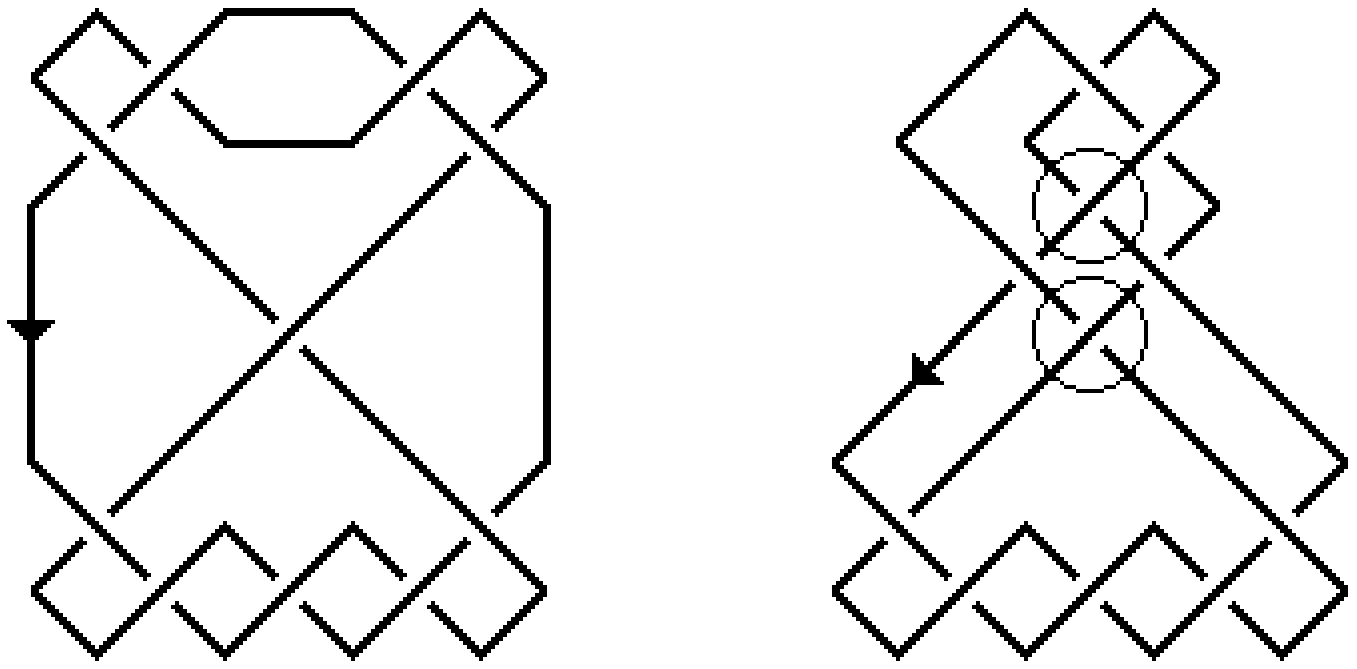

Figure 4 Bleiler-Nalsanishi Unknot Problem

The minimum braid projection has fewer crossings than other projections that have been shown to have a special unknotting number of two for this knot. In his 1994 book, Adams uses a 14 crossing projection with a special unknotting number of two. [16]

When the minimum braid unknotting numbers are compared with the best known values of the classical unknotting numbers [17,18], the numbers are consistent through ten crossing knots. With one exception, whenever the classical unknotting number is only known to be in a range, the minimum braid unknotting number is the upper value of this range. The single exception is 10:1n145 with the minimum braid AAABaaBBCbC. The classical unknotting number for this knot ($10_{131}$) is known to be either one or two. Since only one of the initial A crossings in the minimum braid has to be reversed in order to get an 11 crossing projection of the unknot, the classical unknotting number must be one.

**7. Periodic Tables**

By systematically adding crossings to minimum braids, it is possible to create periodic tables for knots and orientated links. In almost all cases, adding an A crossing to the beginning of a minimum braid creates a minimum braid for another link. This is the case for all knots with ten or fewer crossings and all orientated links with nine or fewer crossings.

The first exception found is a two-component, ten-crossing link with the minimum braid AAbaabcBDcD. Adding an initial A crossing to this braid creates a nonminimum braid for an 11-crossing knot which has a minimum braid AAAbacBaBcb. A link that has a minimum braid that does not generate another minimum braid when an A crossing is added is called a Y-star link.

Because of size limitations, the periodic table of knots and links has been divided into three separate tables and the tables are complete only through seven link crossings. The first table includes knots and two-component links. The second table includes links with two or three components. The third table is a partial listing of "trivial" links. These consist of composite knot and/or links and links with two or more components that are connected by only a pair of crossings. Additional tables could be constructed to include three and four components, four and five components, and so forth.

In each of these tables, adding an A crossing to the beginning of a minimum braid creates a new minimum braid for a different link. This addition is the difference between any two cells in a column of any of the periodic tables. The additional braid crossing almost always changes the number of components by one and increases the number of link crossings by one.

There are two exceptions in the links listed in Tables 1-3. When an A crossing is added to the minimum braid (AbaCbaCbdCd) for link 9:2-41b, the new minimum braid is for a three-component, 11 (not 10) crossing link. The second exception is the minimum braid (AbaCbaCdCbCd) for knot 10:1-113. Adding an initial A crossing to this braid generates a two-component, 12 (not 11) crossing link. A links that have a minimum braid that does not change the number of link crossing by one when an initial A braid crossing is called a Z-star link. Y-star and Z-star links have many interesting properties and will be discussed in another paper.

The first periodic table has a list of numbers on the left-hand side that are the crossing numbers for the links. The data cells are arranged in columns with the two-component links having shaded cells. In each cell there are four rows of data. The top row is the Tag # for the link as listed in Tables 1 and 2.

The second row is an abbreviated minimum braid, where a number in inserted into the minimum braid whenever there are more than two identical crossing in a sequence. For example, the minimum braid for knot 7:1-05 is listed as A4BaBB instead of AAAABaBB.

The third row lists codes for the coefficients of the Alexander polynomial where capital letters represent positive integers and lower case letter represent negative integers. An underscore represents a zero coefficient and is used in the second periodic table. Because the coefficients are symmetric, only half of the coefficients are coded.

Since the Alexander polynomial for links with even (odd) numbers of components have even (odd) numbers of coefficients, two different links can have the same codes. In the first data column, the Alexander polynomial code for knot 3:1-01 and link 4:2-01a is Aa. Since any knot has one component, the coefficients for the trefoil knot are +1, -1, and +1. The coefficients for the two component link 4:2-01a are +1, -1, +1, and -1.

The fourth row is a modified Conway notation [19], where the modification is used to distinguish between links with different orientations and to rearrange the terms in the notation so that the first term always increases by one as one moves down the cells in any column. Whenever the Conway notation has been modified, it is listed in bold type. For example, the third column's first cell is the two-component, four-crossing link 4:2-01b. The other orientation of this link is listed as the third cell in the first column. It has the normal Conway notation of 4, which is consistent with the other Conway notations in the first column. The Conway notations for the third column are m3, where m increase by one in each successive cell.

At the top of each column there is a code that is labeled the column type and consists of five parts. The first part is the number of strands in the minimum braids for that column. The second part is a letter where "a" represents an alternating braid and "n" represents a nonalternating braid. The third part is + (-) if the sums of the Alexander polynomial coefficients for knots in that column are positive (negative) and have a digital equal to 1 (8). The sum of these coefficients is always zero for the links. The fourth part is the z number that is listed in Tables 1 and 2 and is the common exponent of the t variable in the Alexander polynomial. The fifth part is an "o" ("e") if the least significant coefficient of the Alexander polynomial is odd (even).

The three column type parts that are related to Alexander's polynomial are not always constant. The exceptions called Z-star can have a z value that changes between the top two cells (hence the reason they are call Z-stars). There are also some ten crossing links and one 11 crossing knot where the least significant coefficient of the Alexander polynomial can be even or odd in different cells in the same column. There exception are called O-stars and will also be discussed in a future paper.

# Periodic Table of Knots & Links

**1 & 2 Components**

**Column Type**

| Tag #<br>Minimum Braid<br>Alexander Polynomial<br>Conway Notation | |
|---|---|
| **1** | **2** |
| (white) | (shaded) |

| | **2a+0o** | **3a-0o** | **3n+1e** | **3a+0o** | **3n+1e** | **4a-0o** | **4n-1e** | **4a+0o** | **4n-1e** | **4n-1e** | **4n+2o** | **3a+0o** | **4n+2e** | |
|---|---|---|---|---|---|---|---|---|---|---|---|---|---|---|
| **2** | 2:2-01<br>AA<br>A<br>2 | | | | | | | | | | | | | |
| **3** | 3:1-01<br>A3<br>Aa<br>3 | | | | | | | | | | | | | |
| **4** | 4:2-01a<br>A4<br>Aa<br>4 | 4:1-01<br>AbAb<br>Ac<br>**112** | 4:2-01b<br>AABaB<br>B<br>**13** | | | | | | | | | | | |
| **5** | 5:1-01<br>A5<br>AaA<br>5 | 5:2-01<br>AAbAb<br>Ac<br>212 | 5:1-02<br>A3BaB<br>Bc<br>**23** | | | | | | | | | | | |
| **6** | 6:2-01a<br>A6<br>AaA<br>6 | 6:1-01<br>A3bAb<br>AcC<br>312 | 6:2-03<br>A4BaB<br>Bc<br>33 | 6:1-02<br>AAbAbb<br>AcE<br>2112 | 6:2-02b<br>A3BaBB<br>Bd<br>222 | 6:2-02a<br>AbACbC<br>Ae<br>**1122** | 6:1-03<br>AABacBc<br>Be<br>**132** | | | | 6:2-01b<br>AABaBCbC<br>C<br>**15** | | | |
| **7** | 7:1-01<br>A7<br>AaAa<br>7 | 7:2-01a<br>A4bAb<br>AcC<br>412 | 7:1-04<br>A5BaB<br>BcC<br>43 | 7:2-03a<br>A3bAbb<br>AcE<br>3112 | 7:1-05<br>A4BaBB<br>BdE<br>322 | 7:1-02<br>AAbACbC<br>AeG<br>**2122** | 7:2-06<br>A3BacBc<br>Bf<br>232 | 7:1-03<br>AbAbCbC<br>AeI<br>**111112** | 7:2-03b<br>AABaBcBc<br>Bg<br>**12112** | 7:2-01b<br>AAbAbcBc<br>Be<br>**214** | 7:1-06<br>A3BaBCbC<br>Ce<br>**25** | 7:2-04a<br>AAbAAbb<br>AcF<br>**2,21,2** | 7:1-07<br>AABaBBCbC<br>Dg<br>**1213** | |
| **8** | 8:2-01a<br>A8<br>AaAa<br>8 | 8:1-01<br>A5bAb<br>AcCc<br>512 | 8:2-11a<br>A6BaB<br>BcC<br>53 | 8:1-02<br>A4bAbb<br>AcEe<br>4112 | 8:2-12a<br>A5BaBB<br>BdE<br>422 | 8:2-02a<br>A3bACbC<br>AeG<br>3122 | 8:1-10<br>A4BacBc<br>BfG<br>332 | 8:2-04<br>AAbAbCbC<br>AeK<br>211112 | 8:1-12<br>A3BaBcBc<br>BhK<br>22112 | 8:1-11<br>A3bAbcBc<br>BeE<br>**314** | 8:2-11b<br>A4BaBCbC<br>Ce<br>**35** | 8:1-04<br>A3bAAbb<br>AcFg<br>3,21,2 | 8:2-02b<br>A3BaBBCbC<br>Di<br>**2213** | + |
| **9** | 9:1-01<br>A9<br>AaAaA<br>9 | 9:2-01a<br>A6bAb<br>AcCc<br>612 | 9:1-18<br>A7BaB<br>BcCc<br>63 | 9:2-03a<br>A5bAbb<br>AcEe<br>5112 | 9:1-19<br>A6BaBB<br>BdEe<br>522 | 9:1-02<br>A4bACbC<br>AeGg<br>4122 | 9:2-22a<br>A5BacBc<br>BfG<br>432 | 9:1-04<br>A3bAbCbC<br>AeKm<br>311112 | 9:2-24a<br>A4BaBcBc<br>BhL<br>32112 | 9:2-23<br>A4bAbcBc<br>BeE<br>414 | 9:1-29<br>A5BaBCbC<br>CeE<br>**45** | 9:2-05a<br>A4bAAbb<br>AcFg<br>4,21,2 | 9:1-30<br>A4BaBBCbC<br>DiK<br>3213 | + |
| **Al** | 1 : 2 | 3 : 5 | 3 : 4 | 5 : 13 | 5 : 12 | 7 : 12 | 7 : 9 | 13 : 21 | 13 : 18 | 5 : 14 | 5 : 6 | 7 : 20 | 11 : 15 | |
| **Jo** | AAA | a_a_a | A_A_A | a_b_aA | aA_B_A | A_BaAaA | aAaAb_a | aAbBcAbA | aBaCbBaA | A_AaA_A | A_AaA_A | B_BaAa | AaAcAbAa | |
| **Hx** | -1, -2 | -2, -1 | 2, 1 | 1, 2 | 2, 1 | -1, 1 | 1, 2 | 1, 2 | -1, 1 | -1, -2 | -1, -2 | -1, -2 | 1, -1 | |
| **Hr** | 1,0 +1c | -1,0 -0 | 1,0 +2 | 2,1 +3 | 2,1 +4 | 1,0 +1b | -1,0 -1c | 1,0 +1c | -1,0 -1b | -2,1 -3 | 2,1 +3 | 1,0 +5 | 3,2 +5 | |
| **Wy** | 1 | 0 | 3 | 0 | 3 | 1 | 2 | 1 | 2 | -2 | 5 | 0 | 5 | |

At the bottom of Period Table 1, there are five rows that list some interesting properties of the knots and links in each column. The row labeled Al has two numbers that are separated by a colon. The number on the right side is the sum of the absolute values of the coefficients of the Alexander polynomial for the link in the top cell. The number on the left is the difference of this sum between any two cells in a column.

Consider the top cell in the fourth column that is the knot 6:1-02 ($6_3$) and has an Alexander polynomial

$$AP(6{:}1\text{-}02) = t^4 - 3t^3 + 5t^2 - 3t + 1.$$

The minimum braid is AAbAbb; therefore, the number of strands is three and the braid is alternating. The sum of the coefficients of the Alexander polynomial is +1. The z value is zero, and when the polynomial is evaluated for t = 10, the Alexander polynomial equals 10501 - 3030 or 7471. The least significant coefficient is one, so the polynomial is odd and the column type is 3a+Oo. The sum of the absolute values of the coefficients is 13. This sum for the next cell is 18 so the difference between cells is 5, as shown in the row labeled Al.

The next row is labeled Jo and includes a coded description of how the coefficients of the Jones polynomial change between the first and second cell. The Jones polynomial for knot 6:1-02 and link 7:2-03a are

$$J(6{:}1\text{-}02) = -x^6 + 2x^4 - 2x^2 + 3 - 2x^{-2} + 2x^{-4} - x^{-6} \text{ and}$$
$$J(7{:}2\text{-}03a) = -x^9 + 2x^7 - 3x^5 + 3x^3 - 4x + 2x^{-1} - 2x^{-3} + x^{-5}.$$

Add a zero term at the end of the Jones polynomial for the knot to make the number of coefficients equal. The differences between the coefficients of these two Jones polynomial are

0, 0, -1, 0, -2, 0, -1, and 1.

Ignore beginning zero differences and code the remaining differences where A = +1, a = -1, b = -2, and _ = 0 to get a_b_aA. This is the coded difference listed for column four. The signs of these differences alternates between successive pairs of cells so the next coded difference is A_B_Aa.

The next row is labeled Hx and includes the HOMFLYPT polynomials for first two cells when the two variables (x and y) of the HOMFLYPT are set equal to one. The difference equation for the cells starting in the third row of any column is

$$Hx(r) = Hx(r-1) - Hx(r-2),$$

where r represents a row in a given column. This equation and initial conditions generate a cycle of six. For the fourth column, the cycle is 1, 2, 1, -1, -2, and -1. Other columns have this cycle multiplied by a factor of two.

The row labeled Hr is another sequence that uses the HOMFLYPT polynomial. Let the top cell in a column be written as AZ, where Z includes all the braid word except the initial A. The second cell would be AAZ, and the third cell would be AAAZ. Starting with this third cell, the HOMFLYPT polynomial can be calculated by simply taking a weighted sum of the HOMFLYPT polynomials of the first two braids. Instead of expressing the HOMFLYPT polynomial in terms of x and y, the polynomial can be written in terms of the variables r (remove) and s (switch).

H(AAAZ) = rH(AAZ) + sH(AZ),

where H(..) represents the HOMFLYPT polynomial. In terms of x and y, the variables r equals 1/y and s equals -x/y.

If r and s are set equal to one, this equation becomes a second order difference equation with unitary coefficients and can generate Fibonacci or Lucas type sequences. There are three numbers listed for each column. The first two numbers are the initial conditions for the difference equation, and the third number is the value for the top cell's polynomial.

When the sequence is the Fibonacci series (or a multiple or inverse of it), there are three possible starting positions that have a value of one. These three starting positions are labeled with the suffix a, b, or c. Notice that the Fibonacci series is defined to start with 1 and 0 as the initial values so that the first ten terms are 1,0,1,1,2,3,5,8,13, and 21.

In Periodic Table 1, the first column starts with the third position for one and sixth column start with the second position. An example of a column that uses the first position for one in a Fibonacci sequence has link 8:2-09a and knot 9:1-09 as the top two cells. This column would be included in a expanded table 1 that is complete through eight crossing links.

The last row is labeled Wy and is the difference between the topological writhe and the nullification writhe for the links in a column. This variable is called the remaining writhe and will be discussed in the section on a writhe table of knots and links.

Below the periodic table there are dashed and underscore lines to show that certain columns are naturally related. These relationships are important because they determine how the columns in the periodic table are ordered. Whenever columns have links that share a common braid universe, the columns are arranged together. Furthermore, some of these columns may be paired and share some common properties.

Two columns are considered to be *pairs* if one can transform any link into one column into any link of the second column by reversing the minimum braid and then adding A crossings. Consider the links in the second column which have a minimum braid $A^{n}$bAb, where n is greater or equal to one. Reverse this braid and add m+n A crossings to the beginning to get $A^{m+n}a^{n}$BaB or $A^{m}$BaB. If m is greater than one, the braid is the minimum braid for any link in the third column. Column pairs are connected by a dashed line below the two columns.

If two columns share the same braid universe but are not column pairs, the columns are connected with an underscore line. In Periodic Table 1, columns 8-11 all share braid universes and the first two columns are a pair and the last two columns are another pair.

Some of the common characteristics of column pairs can be read directly from the periodic tables. Since they have common braid universes, they must have the same number of strands. Recall that this is the first number in the column type code at the top of each column. The difference between the absolute values of the Alexander polynomial for two cells are the same for each of the columns in a pair. These differences are 3, 5, 7, 13, and 5 for the first five column pairs. The codes for the difference in the coefficients of the Jones polynomials have their order reversed. For some pairs the signs of the codes are reversed. The code for column four is a_b_aA, and the code for column five is aA_B_A. When the columns have alternating knots and orientated links, the sum of the remaining writhes equals three for each pair. These common characteristics of column pairs have not been proved to be general and exceptions might be found as the periodic tables are extended to links with more crossings.

Three other properties of Periodic Table 1 need to be pointed out. There is a space between the first and second columns to isolate column one. This is the only column with two strands and the knots in this column are torus knots. There is an empty column of narrow cells between the last two columns that have data. This column represents column pairs for the last two columns that have been omitted because of lack on space on the page. Finally, there are + signs to the right of the rows for eight and nine crossings to indicate that the table in only complete through seven crossing links.

The second periodic table includes links with two or three components. The cells for the three component links have the top row shaded. Since the sum of the coefficients for Alexander polynomials is equal to zero for all links, the third part of the column type code has been modified. This code now represents the weighted sum of these coefficients and an = sign indicates that this weighted sum is equal to zero. Since this table includes nonalternating links in columns 2-4 and 6, there is an "n" listed for the remaining writhe at the bottom of these columns. Columns 9 and 12 do not have column pairs and the last column has unique braid universes.

The third periodic table includes "trivial" links that are not normally consider as primary knots or links. This table is not complete and only includes one example for each composite or "loop". For example, the top cell in the second column is the square knot. The granny knot is not included since it is another example of the six crossing composite knot. Cells that have four component links have the second row shaded. The last row in each cell is the description of the link as listed in Table 3. Not included at the bottom of the table is a row of remaining writhes and dashed lines to indicate column pairs.

When determining how to order the columns in each of these periodic tables, several guidelines were followed. These guidelines will be described for Periodic Table 1. When building the table, the links with the fewest number of crossing were added first. Since the Hopf link (2:2-01) is the only two-crossing link, it was the first cell added. Adding A crossings to the beginning of the minimum braid for this link determined the remaining cells in the first column. Since the trefoil was included as the second cell in this column,

# Periodic Table of Knots & Links

**2 & 3 Components**

**Column Type**

| Tag #<br>Minimum Braid<br>Alexander Polynomial<br>Conway Notation | |
|---|---|
| **2** | **3** |
| | |

| | 3a=0o | 3n=1o | 3n-1o | 3n-0o | 3a=0o | 3n=1o | 4a=0o | 4n=1e | 4n-1o | 3n-0o | 3n-1e | 3a=0o | |
|---|---|---|---|---|---|---|---|---|---|---|---|---|---|
| **6** | 6:3-01a<br>AAbAAb<br>AcD<br>2,2,2 | | 6:3n03a<br>AAbaab<br>Ab<br>2,2,2- | 6:3n03b<br>AABAAB<br>Aa_<br>2,2,2- | 6:3-02<br>AbAbAb<br>AdF<br>.1 | | | | 6:3-01b<br>AbaCBBC<br>Cf<br>**11,2,2** | | | | |
| **7** | 7:2-02<br>A3bAAb<br>AcD<br>3,2,2 | 7:2n08<br>A3BaaB<br>Ac<br>21,2,2- | 7:2n07a<br>A3baab<br>Aa<br>3,2,2- | 7:2n07b<br>A3BAAB<br>Aa_<br>3,2,2- | 7:2-05<br>AAbAbAb<br>AdG<br>.2 | | 7:3-01a<br>AbACbbC<br>AeH<br>**11,2,2+** | | 7:2-04b<br>AAbACBBC<br>Cg<br>21,2,2 | | 7:3-01b<br>A3BBaBB<br>BeF<br>2,2,2+ | | |
| **8** | 8:3-01a<br>A4bAAb<br>AcDd<br>4,2,2 | 8:3n08a<br>A4BaaB<br>AcD<br>31,2,2- | 8:3n07a<br>A4baab<br>Aa_<br>4,2,2- | 8:3n07b<br>A4BAAB<br>Aa__<br>4,2,2- | 8:3-02a<br>A3bAbAb<br>AdGh<br>.3 | 8:3n09a<br>A3BaBaB<br>AdF<br>(2,2)(2,2-) | 8:2-05<br>AAbACbbC<br>AeJ<br>21,2,2+ | 8:3-06a<br>A3BacBBc<br>BgJ<br>2,2,2++ | 8:3-03b<br>A3bACBBC<br>CgH<br>31,2,2 | 8:3-03a<br>A3bbAbb<br>AcFh<br>31,2,2 | 8:2-09a<br>A4BBaBB<br>BeG<br>3,2,2+ | 8:3-04a<br>AAbAbbAb<br>AdHj<br>.2:20 | + |
| **9** | 9:2-02<br>A5bAAb<br>AcDd<br>5,2,2 | 9:2n44<br>A5BaaB<br>AcD<br>41,2,2- | 9:2n43a<br>A5baab<br>Aa_<br>5,2,2- | 9:2n43b<br>A5BAAB<br>Aa__<br>5,2,2- | 9:2-06<br>A4bAbAb<br>AdGh<br>.4 | 9:2n49<br>A4BaBaB<br>AdG<br>(3,2)(2,2-) | 9:3-01a<br>A3bACbbC<br>AeJl<br>31,2,2+ | 9:2-25<br>A4BacBBc<br>BgK<br>3,2,2++ | 9:2-08b<br>A4bACBBC<br>CgH<br>41,2,2 | 9:2-08a<br>A4bbAbb<br>AcFh<br>41,2,2 | 9:3-09b<br>A5BBaBB<br>BeGh<br>4,2,2+ | 9:2-13<br>A3bAbbAb<br>AdHk<br>.3:20 | + |
| **Al** | 4 : 12 | 4 : 8 | 0 : 4 | 0 : 4 | 8 : 16 | 8 : 16 | 12 : 20 | 12 : 28 | 8 : 12 | 8 : 28 | 8 : 20 | 12 : 36 | |
| **Jo** | bab_a | A_BAB | aaaaaa | aaaaaa | b_c_bA | Ab_c_b | a_cAcAbA | AbAcAc_a | a_c_bAa | a_c_bAa | aAb_c_a | b_dAbBa | |
| **Hx** | -2, -1 | -1, -2 | 1, -1 | -2, -1 | 1, 2 | 1, 2 | -2, -1 | 1, -1 | 1, -1 | 1, 2 | 1, -1 | 4, 2 | |
| **Hr** | -1,0 -2 | -1,0 -1c | -2,-1 -1 | 2,1 +4 | 1,0 +1c | -1,0 -1c | 1,0 +0 | -1,0 -3 | 1,0 +3 | -3,-2 -7 | 3,1 +5 | -2,0 -4 | |
| **Wy** | 0 | n | n | n | 0 | n | 1 | 2 | 4 | 0 | 3 | 0 | |

# Periodic Table of Loops & Composites

**Column Type**

| Tag #<br>Minimum Braid<br>Alexander Polynomial<br>Factor Knots / Links | | | |
|---|---|---|---|
| 1 | 2 | 3 | 4 |

| | 3a-0o | 3a+0o | 4a=0o | 4a=0o | 4a-0o | 4a+0o | 4n-1e | 4n-1e | |
|---|---|---|---|---|---|---|---|---|---|
| 4 | 4:3.01<br>AAbb<br>Ab<br>L2+1 | | | | | | | | |
| 5 | 5:2.01<br>A3bb<br>Ab<br>K3+1 | | | | | | | | |
| 6 | 6:3.01a<br>A4bb<br>AbB<br>L4a+1 | 6:1c01a<br>A3b3<br>AbC<br>K3#K3 | 6:4.01<br>AAbbCC<br>Ac<br>L2+1+1 | 6:4.02<br>AAbCCb<br>Ac<br>1+L2+1 | 6:2.01<br>AAbCbC<br>Ad<br>1+K4 | | 6:3.01b<br>AABaBCC<br>Bd<br>L4b+1 | | |
| 7 | 7:2.01<br>A5bb<br>AbB<br>K5.1+1 | 7:2c01<br>A4b3<br>AbC<br>L4a#K3 | 7:3.01<br>A3bbCC<br>AcD<br>K3+1+1 | 7:3.02<br>A3bCCb<br>AcD<br>1+K3+1 | 7:1c01<br>A3bCbC<br>AdE<br>K3#K4 | 7:3.03<br>AAbAbCC<br>AdF<br>L5+1 | 7:2.02<br>A3BaBCC<br>Be<br>K5.2+1 | 7:2c02<br>A3bbcBc<br>Bd<br>K3#L4b | |
| 8 | 8:3.01a<br>A6bb<br>AbBb<br>L6.1a+1 | 8:1c01<br>A5b3<br>AbCc<br>K5.1#K3 | 8:4.01a<br>A4bbCC<br>AcD<br>L4a+1+1 | 8:4.02a<br>A4bCCb<br>AcD<br>1+L4a+1 | 8:2c01<br>A4bCbC<br>AdE<br>L4a#K4 | 8:2.01<br>A3bAbCC<br>AdF<br>K6.1+1 | 8:3.04<br>A4BaBCC<br>BeF<br>L6.3+1 | 8:3c02<br>A4bbcBc<br>BdD<br>L4a#L4b | + |
| 9 | 9:2.01<br>A7bb<br>AbBb<br>K7.1+1 | 9:2c01<br>A6b3<br>AbCc<br>L6.1a#K3 | 9:3.01<br>A5bbCC<br>AcDd<br>K5.1+1+1 | 9:3.02<br>A5bCCb<br>AcDd<br>1+K5.1+1 | 9:1c01<br>A5bCbC<br>AdEe<br>K5.1#K4 | 9:3.03<br>A4bAbCC<br>AdFf<br>L7.1a+1 | 9:2.02<br>A5BaBCC<br>BeF<br>K7.4+1 | 9:2c02<br>A5bbcBc<br>BdD<br>L5#L4b | + |
| **Al** | 2 : 4 | 3 : 9 | 4 : 8 | 4 : 8 | 5 : 10 | 6 : 16 | 6 : 8 | 4 : 12 | |
| **Jo** | aabaa | aab__A | AACBCAA | AACBCAA | A_AaA_A | a_b_b_a | a_b_b_a | A_A_BAA | |
| **Hx** | 1, 2 | 4, 2 | -1, -2 | -1, -2 | 2, 4 | 1, -1 | -2, -1 | -4, -2 | |
| **Hr** | -1,0 -1c | 2,0 +4 | -1,0 -1c | -1,0 -1c | 0,0 0 | -1,0 -1b | 1,0 +2 | -2,0 -4 | |

the next two columns added start with a four crossing knot (4:1-01) or orientated link (4:2-01b). Since these two columns are a column pair and have the same number the strands, the column with alternating minimum braids was included first.

The next decision was how to arrange the columns that start with six crossing links that had not already been included in one of the first three completed columns. The five links are 6:1-02, 6:2-02b, 6:2-02a, 6:1-03, and 6:2-01b. In general, when determining the ordering for a given number of link crossing, the links with the fewest number of strands are included first. Column pairs are always included together along with columns that are not pairs but have common braid universes. For column pairs, alternating braids are always included before nonalternating braids. Following these guidelines, the next five columns in the first period table are added in the order that the six crossing links were listed previously.

When extending this periodic table to seven crossing links, there are five links that had not already been included in of the first eight completed columns. Three of these columns have common braid universes with the column headed by 6:2-01b. Since the column headed by knot 7:1-03 has alternating minimum braids, this column and its column pair (7:2-03b) are included *before* the column for 6:2-01b. Since columns headed by 7:2-01b and 6:2-01b are a pair, it is necessary to determine which one should be included first. The column headed by 7:2-01b was selected because it has the longest initial sequence of alternating braid crossings. This ordering appears to be consistent with the fourth criterion for determining minimum braids. The final two columns added are 7:2-04a and 7:1-07 with the column headed by 7:2-04a included first since it has three strands.

**8. Applications of Periodic Tables.**

Besides providing a way to organize the data that are included in them, the periodic tables can be applied to study other properties of knots and orientated links. Three such applications will be described in the following section.

The first application studies how to draw some of the links so as to show how they are related. Consider the links that are included in the sixth column of Periodic Table 1. The top line of Figure 5 plots these links as they were drawn by Roth in Rolfsen's Knots and Links [3]. From these drawings it is not very clear how these links are related to each other. The bottom line of Figure 5 plots these same links using their minimum braids. In each case the right part of the drawing remains the same. As additional A crossings are added to the minimum braids, the number of crossings in the left part of the drawing increases by one.

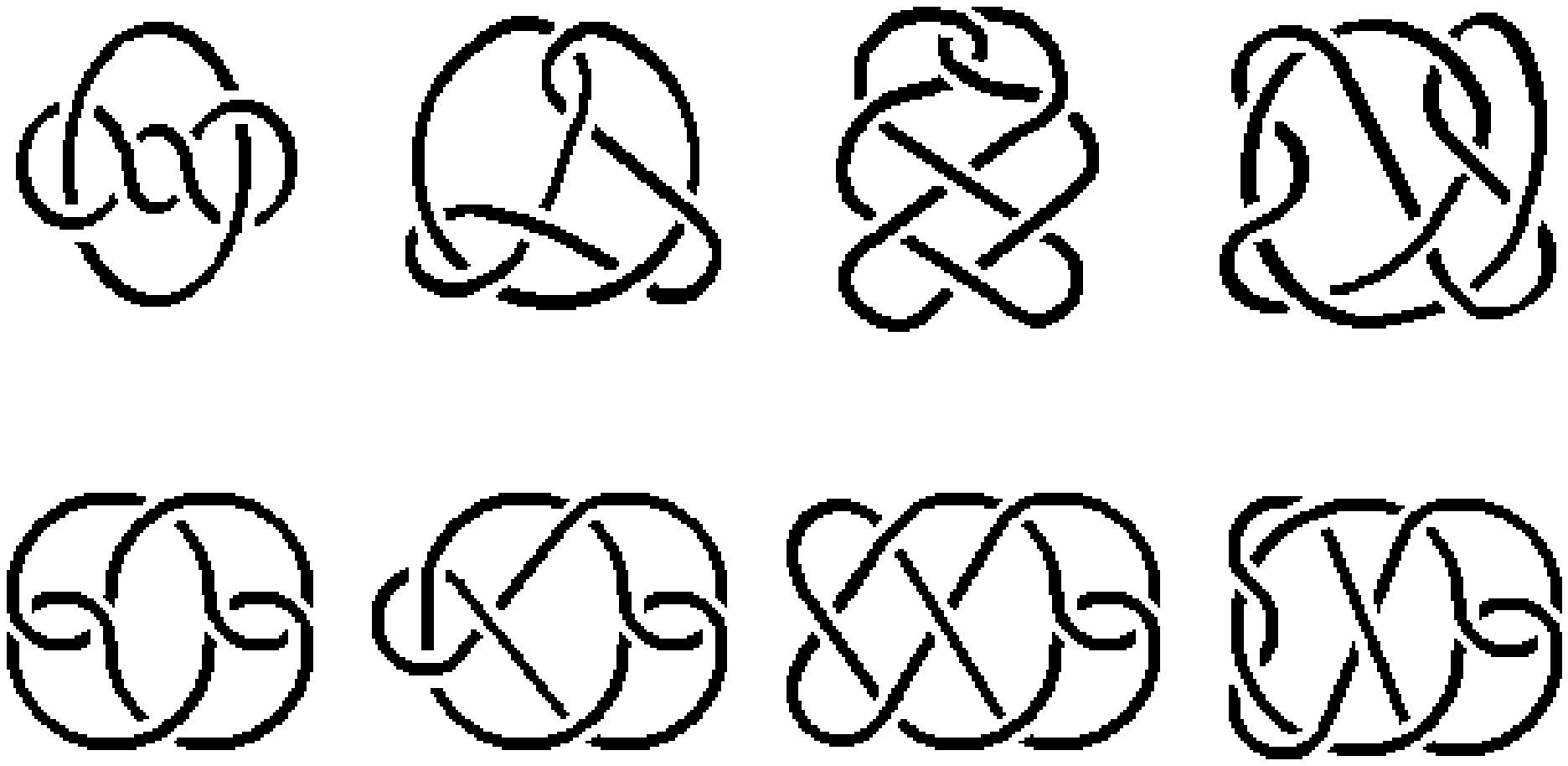

Figure 5. Links of the Sixth Column

A second application involves the study of hyperbolic volumes of links. Since adding A crossings to the beginning of a minimum braid to create columns in a periodic table is a simple type of Dehn surgery [20], the periodic tables can be used to study hyperbolic volumes. Figure 6 is a graph of the hyperbolic volumes of links in the first three column pairs of Periodic Table 1. This graph is divided into two parts. The left section plots the hyperbolic volumes for the first five links in columns 2-7. M is used to specify the link position and is equal to one when the link is in the top cell of a column. The hyperbolic volumes were calculated in the SnapPea program by Weeks. The hyperbolic volumes for knots and links are listed in [21].

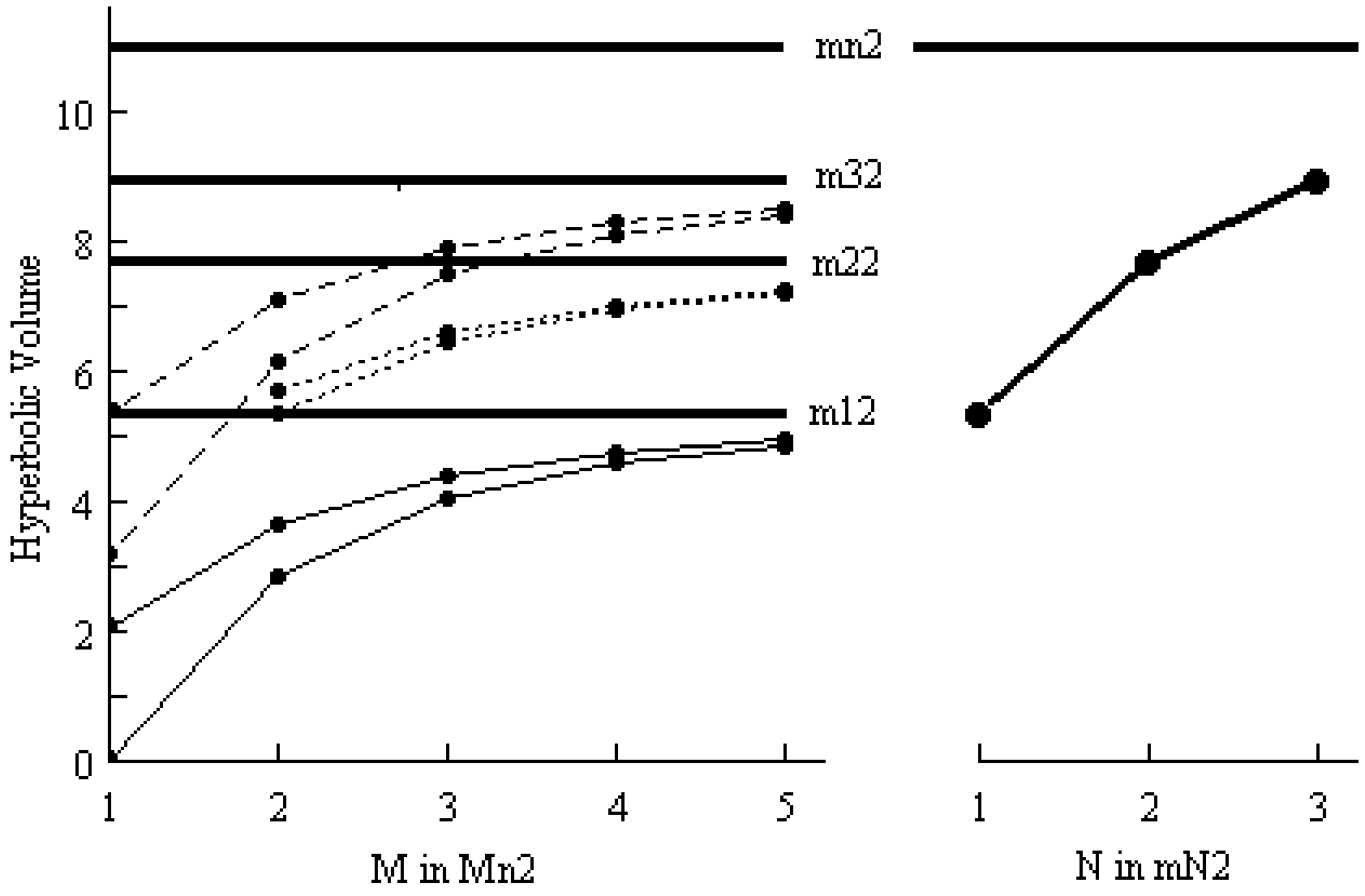

Figure 6. Asymptotic Hyperbolic Volumes

The hyperbolic volumes for the links in the second and third columns are plotted as solid lines that are asymptotically approaching the horizontal line labeled m12. The fourth and fifth columns are plotted as dotted lines and their asymptote is m22. The sixth and seventh columns are plotted as dashed lines and their asymptote is m32. The labels for these three asymptotes are simply the Conway notations for columns 2, 5, and 7. On the right side of Figure 6 the asymptotic values are plotted where N represents the second integer of the Conway notation for these rational links. The horizontal line labeled mn2 is the asymptote for the asymptotic values as N increases to infinity. Of course, this procedure could be extended to include the asymptotic values for mn3, mn4, and so on. The "triple" asymptote mno would be the largest hyperbolic volume of rational links with three integer Conway notations.

A third application of periodic tables is the study of ideal writhes. [22,23,24] Figure 7 plots the ideal writhe of links with minimal crossing numbers (MCN) of seven or less. Since knots and orientated links are plotted, the ideal writhes are expressed in units of 2/7. In this plot there are two solid lines drawn with a slope equal to five. These lines represent moving down a column in Periodic Table 1. The solid line that has a MCN of two is the Hopf link in cell one of the first column. The solid line that has a MCN of four and an ideal writhe of zero represents moving down the second column in Periodic Table 1.

The two dashed lines have slopes of two and connect the twist knots. The top line connects the twist knots with odd numbers of crossings and the bottom line connects the twist knots with even numbers of crossings. The dotted line that also has a slope of two connects the "ABC" oriented links. The abbreviated Alexander polynomials of these links are A, B, and C, respectively.

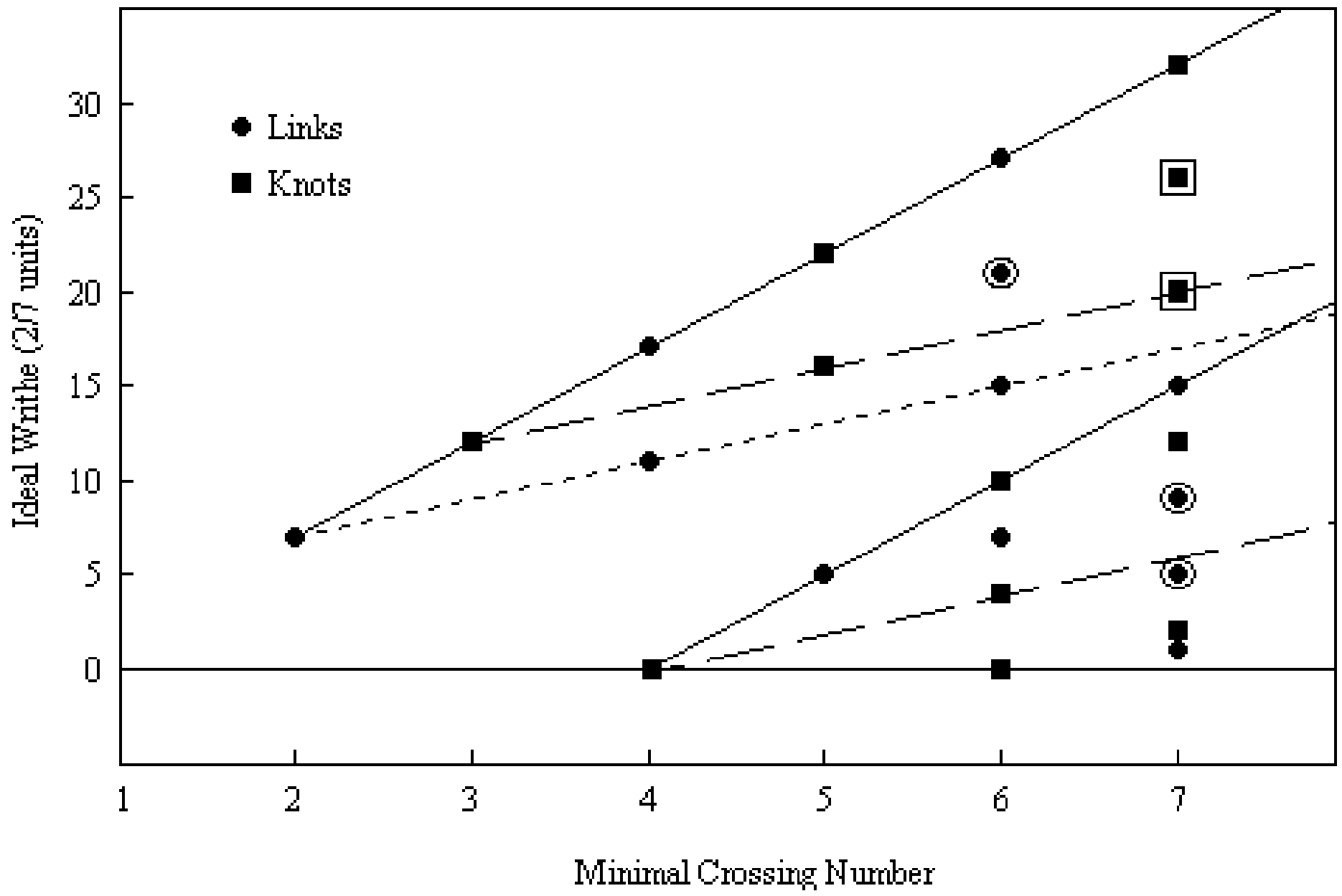

Figure 7. Ideal Writhes of Links with 2-7 Minimal Crossing Numbers

This plot of ideal writhes can be used to reorder the columns to create a Writhe Table of Knots and Links. The new order depends upon the position of the ideal writhes for a minimum crossing number of seven. Start with the top most link in Figure 7 and move down to determine the left-to-right ordering in the new writhe table. When there are pairs of knots or links with the same ideal writhes and a MCN of seven, the two links are assigned the same relative order as in Periodic Table 1.

The writhe table is organized like the periodic tables with some differences to be noted. The second row for each cell lists the link's Thurston-Bennequin numbers separated by a colon. These numbers were calculated from the Kauffman polynomials [25] for the link and its reverse using Ochiai's KnotTheorybyComputer program for the Macintosh. The third row in each cell lists the topological writhe and the nullification writhe, separated by a colon. The writhes are consistent with the definition of a positive crossing that is used throughout this paper. The fourth row is the ideal writhe measured in units of 2/7. With these units, the ideal writhe is equal to two times the topological writhe plus three times the nullification writhe.

At the bottom of the table is a row that lists the remaining writhe. By definition, the remaining writhe is the difference between the topological writhe and the nullification writhe. This row is labeled Wy and is included at the bottom of the periodic table.

The columns are divided into four sections that correspond to the plot of the ideal writhe. All the links that are between the upper solid line and the dotted line in Figure 7 are included in the first group. The links that are be between the lower solid and dashed lines are in the second group.

# Writhe Table of Knots & Links

**1 & 2 Components**

| Tag #<br>Thurston-Bennequin #'s<br>Writhe (-) Reduction Writhe<br>Ideal Writhe (2/7 units) |
|---|

| | 2a+0o | 3n+1e | 3n+1e | 4n+2o | 4n+2e | 3a-0o | 4a-0o | 4n-1e | 4n-1e | 3a+0o | 3a+0o | 4a+0o | 4n-1e | |
|---|---|---|---|---|---|---|---|---|---|---|---|---|---|---|
| **2** | 2:2-01<br>0 : -4<br>2 : 1<br>7 | | | | | | | | | | | | | |
| **3** | 3:1-01<br>1 : -6<br>3 : 2<br>12 | | | | | | | | | | | | | |
| **4** | 4:2-01a<br>2 : -8<br>4 : 3<br>17 | 4:2-01b<br>0 : -6<br>4 : 1<br>11 | | | | 4:1-01<br>-3 : -3<br>0 : 0<br>0 | | | | | | | | |
| **5** | 5:1-01<br>3 : -10<br>5 : 4<br>22 | 5:1-02<br>1 : -8<br>5 : 2<br>16 | | | | 5:2-01<br>-2 : -5<br>1 : 1<br>5 | | | | | | | | |
| **6** | 6:2-01a<br>4 : -12<br>6 : 5<br>27 | 6:2-03<br>2 : -10<br>6 : 3<br>21 | 6:2-02b<br>2 : -10<br>6 : 3<br>21 | 6:2-01b<br>0 : -8<br>6 : 1<br>15 | | 6:1-01<br>-1 : -7<br>2 : 2<br>10 | 6:2-02a<br>-2 : -6<br>2 : 1<br>7 | 6:1-03<br>-3 : -5<br>2 : 0<br>4 | | 6:1-02<br>-4 : -4<br>0 : 0<br>0 | | | | |
| **7** | 7:1-01<br>5 : -14<br>7 : 6<br>32 | 7:1-04<br>3 : -12<br>7 : 4<br>26 | 7:1-05<br>3 : -12<br>7 : 4<br>26 | 7:1-06<br>1 : -10<br>7 : 2<br>20 | 7:1-07<br>1 : -10<br>7 : 2<br>20 | 7:2-01a<br>0 : -9<br>3 : 3<br>15 | 7:1-02<br>-1 : -8<br>3 : 2<br>12 | 7:2-06<br>-2 : -7<br>3 : 1<br>9 | 7:2-03b<br>-2 : -7<br>3 : 1<br>9 | 7:2-03a<br>-3 : -6<br>1 : 1<br>5 | 7:2-04a<br>-3 : -6<br>1 : 1<br>5 | 7:1-03<br>-4 : -5<br>1 : 0<br>2 | 7:2-01b<br>-4 : -5<br>-1 : 1<br>1 | |
| **8** | 8:2-01a<br>6 : -16<br>8 : 7<br>37 | 8:2-11a<br>4 : -14<br>8 : 5<br>31 | 8:2-12a<br>4 : -14<br>8 : 5<br>31 | 8:2-11b<br>2 : -12<br>8 : 3<br>25 | 8:2-02b<br>2 : -12<br>8 : 3<br>25 | 8:1-01<br>1 : -11<br>4 : 4<br>20 | 8:2-02a<br>0 : -10<br>4 : 3<br>17 | 8:1-10<br>-1 : -9<br>4 : 2<br>14 | 8:1-12<br>-1 : -9<br>4 : 2<br>14 | 8:1-02<br>-2 : -8<br>2 : 2<br>10 | 8:1-04<br>-2 : -8<br>2 : 2<br>10 | 8:2-04<br>-3 : -7<br>2 : 1<br>7 | 8:1-11<br>-3 : -7<br>0 : 2<br>6 | + |
| **9** | 9:1-01<br>7 : -18<br>9 : 8<br>42 | 9:1-18<br>5 : -16<br>9 : 6<br>36 | 9:1-19<br>5 : -16<br>9 : 6<br>36 | 9:1-29<br>3 : -14<br>9 : 4<br>30 | 9:1-30<br>3 : -14<br>9 : 4<br>30 | 9:2-01a<br>2 : -13<br>5 : 5<br>25 | 9:1-02<br>1 : -12<br>5 : 4<br>22 | 9:2-22a<br>0 : -11<br>5 : 3<br>19 | 9:2-24a<br>0 : -11<br>5 : 3<br>19 | 9:2-03a<br>-1 : -10<br>3 : 3<br>15 | 9:2-05a<br>-1 : -10<br>3 : 3<br>15 | 9:1-04<br>-2 : -9<br>3 : 2<br>12 | 9:2-23<br>-2 : -9<br>1 : 3<br>11 | + |
| **Wy** | 1 | 3 | 3 | 5 | 5 | 0 | 1 | 2 | 2 | 0 | 0 | 1 | -2 | |

There are a number of observations that can be made about the data presented in the writhe table. Since ideal writhes are not defined for nonalternating links and since there are known exceptions (Z-star, Y-star, and O-star) to the periodic tables, these observations are restricted to this particulate writhe table. In the four groups the difference between the topological writhes and link crossing numbers are 0, 4, 6 and 8, respectively. In each cell the absolute value of the sum of the Thurston-Bennequin numbers is equal to the link crossing numbers plus two. The difference between the Thurston-Bennequin numbers is equal to three times the topological writhe minus two times the remaining writhe. These types of regularities help to justify using minimum braids to generate consistent orientations for different knots and orientated links.

## 9. Conclusions

This paper has presented a new invariant for knots and oriented links and has described some applications. This invariant is complete and appears to be very useful in generating families of knots that have some common property. Additional applications that are being studied include strand passage metrics, Morton-Franks-Williams knots, identical HOMFLYPT and other polynomials, and pretzel knots. In two forthcoming papers, minimum braids are generated for knots through 11 crossings and orientated links through 10 crossings [26] and minimum braids are calculated from certain nonminimum braids using Reidemeister-Markov moves. [10]

I hope these papers and tables help somebody find a top-down way to identify minimum braids and calculate them directly from any link projection or braid. This appears to be a very difficult task, but one that is well worth trying to solve.